\documentclass[twoside,12pt]{article}
\usepackage[latin1]{inputenc}
\usepackage{amssymb}
\newtheorem{theorem}{Theorem}
\newtheorem{myth}{Theorem}[section]
\newtheorem{mylem}{Lemma}[section]
\newtheorem{mypro}{Proposition}[section]

\newtheorem{acknowledgement}[theorem]{Acknowledgement}

\begin{document}
\title{\vspace{-1in}\parbox{\linewidth}{\footnotesize\noindent}
%\vspace{\bigskipamount} \\
{\sc Nonlinear Boundary Stabilization for Timoshenko Beam  System}\\{\tiny In memorian ao Prof. Silvano B. Menezes}}
\date{}
\author{\textsc{A. J. R. Feitosa \,\thanks{%
Universidade Federal da Para\'iba,\, CCEN-DM, PB, Brasil,
\,joaquim@mat.ufpb.br} ,\, M. L. Oliveira \,\thanks{%
Universidade Federal da Para\'iba,\, UFPB, DM, PB, Brasil,
\,,milton@mat.ufpb.br}\, M. Milla Miranda \,\thanks{%
Universidade Estadual da Para\'iba,\, DM, Brasil,\,milla@im.ufrj.br }}}
\maketitle
\begin{quote}
{\scriptsize
{\bf abstract:}

This paper is concerned with the existence and decay of solutions of the following Timoshenko system:
$$
\left\|\begin{array}{cc}
u''-\mu(t)\Delta u+\alpha_1 \displaystyle\sum_{i=1}^{n}\frac{\partial v}{\partial x_{i}}=0,\, \in \Omega\times (0, \infty),\\
v''-\Delta v-\alpha_2 \displaystyle\sum_{i=1}^{n}\frac{\partial u}{\partial x_{i}}=0, \, \in \Omega\times (0, \infty),
\end{array}
\right.
$$
subject to the nonlinear boundary conditions,
$$
\left\|\begin{array}{cc}
u=v=0 \,\, in \,\Gamma_{0}\times (0, \infty),\\
\frac{\partial u}{\partial \nu} + h_{1}(x,u')=0\, in\,\, \Gamma_{1}\times (0, \infty),\\
\frac{\partial v}{\partial \nu} + h_{2}(x,v')+\sigma (x)u=0 \, in\, \,\Gamma_{1}\times (0, \infty),
\end{array}
\right.
$$
and the respective initial conditions at $t=0$. Here $\Omega$ is a bounded open set of $\mathbb{R}^n$ with boundary $\Gamma$ constituted by two disjoint parts $\Gamma_{0}$ and $\Gamma_{1}$ and $\nu(x)$ denotes the exterior unit normal vector at $x\in \Gamma_{1}$. The functions $h_{i}(x,s),\,\, (i=1,2)$ are continuous and strongly monotone in $s\in \mathbb{R}$.

The existence of solutions of the above problem is obtained by applying the Galerkin method with a special basis, the compactness method and a result of approximation of continuous functions by Lipschitz continuous functions due to Strauss. The exponential decay of energy follows by using appropriate Lyapunov functional and the multiplier method.}\\[10pt]
{\scriptsize

}
\end{quote}
\section{introduction}
\setcounter{equation}{0}
The small vibrations of an elastic beam of length $L$ when are considered the rotatory inertial and sheared force can be studied by the following system of equations,
\begin{equation}
\begin{array}{cc}
\displaystyle \frac{\partial^2 u}{\partial t^2}(x,t)-cd\frac{\partial^2 u}{\partial x^2}(x,t) -cdv(x,t) = 0,\,\,0<x<L,\; t\geq 0\label{a}\\
\displaystyle \frac{\partial^2 v}{\partial t^2}(x,t)-c\frac{\partial^2 v}{\partial x^2}(x,t) -c^2d\frac{\partial u}{\partial x}(x,t)+c^2dv(x,t) = 0,\,\,0<x<L,\; t\geq 0
\end{array}
\end{equation}
Completed with the boundary conditions
\begin{equation}
\displaystyle u(0,t)=0, v(0,t)=0, \frac{\partial u}{\partial x}(L,t)=0, \frac{\partial v}{\partial x}(L,t)=0, t>0,\, \label{b}
\end{equation}
and initial conditions
\begin{equation}
\displaystyle u(x,0)=u^0, v(x,0)=v^0,\, \frac{\partial u}{\partial t}(x,0)=u^1(x),\frac{\partial v}{\partial t}(x,0)=v^1(x),\,\, 0<x<L \label{c}
\end{equation}
Here $u(x,t)$ and $v(x,t)$ denote the transversal displacement and the rotation, respectively, of the point $x$ of the beam at the instant $t$. In (\ref{a}), $c$ and $d$ represent the constants:
\begin{eqnarray*}
c=\frac{AL^2}{I_1}\,\, \mbox{and} \,\, d=\frac{GI_1}{EI}
\end{eqnarray*}
where $A$ is the cross sectional area, $G$ is the modulus of elasticity and  $E$ is the shear Young modulus, respectively, of the beam.  $I,\,I_1$ are the axial inertial moment and polar moment, respectively.

System (\ref{a}) was introduced by Timoshenko \cite{Timoshenko}. In Tuscnak \cite{Tucsnak} can be found a nonlinear version of (\ref{a}). The boundary condition (\ref{b}) denote that the end $x=0$ of the beam remains fixed and the end  $x=L$, built-in, with the boundary conditions
\begin{equation}
\left\|\begin{array}{cc}
u(0,l)=v(0,l)=0,\,\, t>0;\\
cd\Bigl[ \frac{\partial u}{\partial x}(L,t)-v(L,t)\Bigr]=-\delta\frac{\partial u}{\partial x}(L,t), \,\, t > 0,\,(\delta \,\, constant);\label{d}\\
\frac{\partial v}{\partial x}(L,t) = -\tau \frac{\partial v}{\partial t}(L,t), \, \, t>0,\, (\tau>0\,\,constant).
\end{array}
\right.
\end{equation}
Kim and Renardy \cite{Kim} studied the existence of solutions of (\ref{a}). Tucsnak \cite{Tucsnak} obtained the existence and exponential decay of solutions for this nonlinear version of (\ref{a})  but with small initial data.

Let $\Omega$ be a bounded open set of $\mathbb{R}^n$ with boundary $\Gamma$ constituted by two disjoint parts $\Gamma_{0}$ and  $\Gamma_{1}$, $\overline {\Gamma_{0}}$$\bigcap$$\overline {\Gamma_{1}}$$= \Phi$. By $\nu(x)$ we represent the exterior unit normal vector at $x\in \Gamma_{1}$. A significant generalization of Problem (\ref{a}), (\ref{c}), (\ref{d}) is the following:
\begin{equation}
\left\|\begin{array}{cc}
\displaystyle u^{''}(x,t)-\mu(t)\Delta
u(x,t)+\alpha_1\sum_{i=1}^{n}\frac{\partial v}{\partial x_i}=0,\;x\in\Omega,\, t>0\\
\label{e}
\displaystyle v{''}(x,t)-\Delta v(x,t) -\alpha_2\sum_{i=1}^{n}
 \frac{\partial u}{\partial x_i}=0, \; x\in \Omega \,t>0
\end{array}
\right.
\end{equation}
\begin{equation}
\left\|\begin{array}{cc}
u(x,t)=0,\;\;v(x,t)=0\;x\in \Gamma_{0},\, t>0;\\
\displaystyle\frac{\partial u}{\partial \nu }(x,t)+h_1(x,u'(x,t))=0,\;x\in\Gamma_{1}\, t>0; \label{f}\\
\displaystyle\frac{\partial v}{\partial \nu }(x,t)+h_2(,x,v'(x,t))+ \sigma(x)u(x,t)=0,\;x\in\Gamma_{1},\,t>0;
\end{array}
\right.
\end{equation}

\begin{equation}
u(x,0)=u^0(x),\;v (x,0)=v^0(x),\, u'(x,0)=u^1(x),\; v'(x,0)=v^1(x)\,x\in \Omega\; \label{g}.
\end{equation}

Here $\mu(t),\;\sigma(x),\; h_1(x,u'(x,s)),\; h_2(x,u'(x,s))$  are real functions defined in $t>0,\;x\in \Gamma_{1} $ and $ x\in\mathbb{R}$, respectively, and $\alpha_1,\;\alpha_2$ are constants.

In Mota \cite{Mota} was analyzed the existence and exponential decay of solutions of Problem (\ref{e}) - (\ref{g}). In this work, the author consider a nonlinear version of (\ref{e}) but the boundary conditions on $\Gamma_{1}$ are linear, i. e.; $h_1(x,s)=\delta_1(x)s$ and $h_2(x,s)=\delta_2(x)s$. Of course, the initial data are small.

In the case of wave equation (i. e., when $\mu =1$ and $\alpha =0$ in $(\ref{e})_1$) with linear boundary dissipation on $\Gamma_{1}$ (i. e.; $h(x,s)=\delta(x)s$), Komornik and Zuazua \cite{Komornikzua}, using the semigroup theory, showed the existence of solutions. Under the same conditions, but applying the Galerking method with a special basis, Milla Miranda and Medeiros \cite{Milla} obtained similar results. The second method, furthermore to be constructive, has the advantage of showing the Sobolev space where lies $\displaystyle \frac{\partial u}{\partial\nu}$.

The above second method has been applied with success to obtain existence of solutions of divers equations, first, with linear boundary dissipations and then, for nonlinear boundary dissipations. In the first case, we can mention the papers of Clark et al. \cite{Clark}, for a coupled system; Araruna and Maciel \cite{Araruna}, for the Kirchhoff equation; Mota \cite{Mota}, for a nonlinear Timoshenko system and Araujo et al. \cite{Araujo}, for a beam equation. In the second case, we cite, among others, the works of Louredo and Milla Miranda \cite{LM},for a coupled system of Klein-Gordon equations; Louredo and Milla Miranda\cite{Louredo}, for a coupled system of Kirchhoff equations and Louredo et al. \cite{millaaraujo}, for a nonlinear wave equation.

The existence of solutions of the wave equations with a nonlinear boundary dissipations has been obtained, among other, applying the theory of monotone operators by Zuazua \cite{Zuazua}, Lasiecka and Tataru \cite{Laz}  and Komornik \cite{Komornik}, and applying the Galerkin method by Vitillaro \cite{Vitilaro} and Cavalcanti et al. \cite{Cavalcanti}

In all of the above works, the exponential decay of solutions is obtained by applying a Lyapunov functional and the technique of multipliers, see Komornik and Zuazua \cite{Komornikzua}.

It is worth emphasizing that the known results in the exponential decay of solution  of the wave equation with nonlinear boundary dissipation where obtained by supposing that $h(s)$ has a linear behavior in the infinite, i.e.,
\begin{equation}
d_0|s|\leq |h(s)|\leq d_1|s|\; \forall s\geq R  \label{h}
\end{equation}
$R$ sufficiently large ($d_0$ and $d_1$ positive constants), see Komornik \cite{Komornik} and the references therein.

In this paper we study the existence and decay of solutions of Problem (\ref{e}) - (\ref{g}). In the existence of solutions we consider two general functions $h_i(x,s)\; (i=1,2)$ which are continuous and strongly monotone in s, i. e.,
$$\bigl[h_i(x,s)-h_i(x,r)\bigr] \geq d_i(s-r)^2,\;\forall s,r\in \mathbb{R},\; x\in\Gamma_{1}, (i=1,\;2) $$
In this part we apply the Galerkin method with a special basis, the compactness method and a result of approximation of continuous functions by Lipschitz continuous functions (see Straus [18]). The choice of the special basis allows us to bound the approximate solutions ($u_{lm}$), ($v_{lm}$) of Problem (\ref{e}) - (\ref{g}) at $t=0$. This in turn permits us to pass to the limit in the nonlinear parts ($h_i(.,u'_{lm})$), ($i=1,\;2$). The exponential decay of energy is obtained for particular $h_i(x,s)=\bigl[m(x)\nu(x)\bigr]p_i(s),\; (i=1,\;2)$, where $m(x)=x-x_0$ and $p_i(x)$ is continuous, strongly monotone and satisfies (\ref{h}). In this part we use an appropriate Lyapunov functional and the multiplier method.
It is important to emphasize that initially we do not know if the sign of the derivative of the energy $E (t)$ associated to our system is negative, to overcome this difficulty, we add to it an appropriate functional $F (t)$, so that the derivative of $\frac{d}{dt} (E + F)$ becomes negative and thus we prove that the energy of the studied system decays at an exponential rate. Until now we do not know any work where the sign of derivative of the energy of the system is not known. This is a novelty in our work.

\section{Notations and Main Results}
\setcounter{equation}{0}
Let $\Omega$ be a bounded open set of $\mathbb{R}^n$ with a $C^2$-boundary $\Gamma$ constitutedbe two disjoint parts $\Gamma_0,\;\Gamma_1$ with $\overline {\Gamma_{0}}$$\bigcap$$\overline {\Gamma_{1}}$$= \Phi$ and $mes(\Gamma_0)>0,\;mes(\Gamma_1)>0$. The scalar product and norm of the real Hilbert space  $L^2(\Omega)$ are denoted by $(u,v)$ and $|u|$, respectively. By $V$ is represent  tel the Hilbert space.
$$V=\{u\in H^1(\Omega)  ;\; u=0\; in \; \Gamma_0 \} $$
%%%%%%%%%%%%%%%%%%%%%%%%%%%%%%%%%%%%%%%%%%%%%%%%%%
provided with the  sorts product and norm
%%%%%%%%%%%%%%%%%%%%%%%%%%%%%%%%%%%%%%%%%%%%%%%%%%%
$$((u,\,v )) = \displaystyle\sum_{i=1}^{n}(\frac{\partial u}{\partial x_i},\,\frac{\partial v}{\partial x_i}),\;\; ||u||=((u,\,u))^2  $$

Let $A=-\Delta $ be the self-adjoint operator determined by the triplet $\{V,\,L^2(\Omega), ((,\,)) \}$ (see Lions \cite{LJL}). Then
$$ D(-\Delta)=\{u\in V\cap H^2(\Omega); \, \frac{\partial u}{\partial \nu}=0\,\, on\,\,  \Gamma_1  \}$$

In order to state the result on the existence of solutions, we introduce the necessary hypotheses. Consider functions
\begin{equation}
h_i \in C^0(\mathbb{R}, L^{\infty}(\Gamma_1)),\;\;h_i (x,0)=0\,\, \mbox{a.\,e.},\;x\in\Gamma_1  \label{i}
\end{equation}
$(i=1,\;2)$ which are strongly monotone in the second variable, i.e.,
 \begin{equation}
\bigl[h_i(x,s)-h_i(x,r)\bigr] \geq d_i(s-r)^2,\;\forall s,r\in \mathbb{R} \label{j}
\end{equation}
a.e. $x\in \Gamma_{1}$  where $d_i$ are positive constants $(i=1,\;2)$. Also consider
\begin{equation}
\mu\in W^{1,\,1}_{loc}(0,\,\infty),\;\;\mu(t)\geq \nu_0> 0,\;\;\forall t\in [0,\infty),\; (\nu_0\;\mbox{constant}) \label{k}
\end{equation}
%%%%%%%%%%%%%%%%%%%%%%%%%%%%%%%%%%%%%%%%%%%%%%%%%%%%%%%%%%%%%%%%%%%%%%%%%%%%%%%%%%%%%%%%%%%%%%%%%%%%%%%%%%%%%%%%%%%%%%%%%%%
and
\begin{equation}
\sigma\in W^{1,\;\infty}(\Gamma_{1}) \label{l}
\end{equation}
%%%%%%%%%%%%%%%%%%%%%%%%%%%%%%%%%%%%%%%%%%%%%%%%%%%%%%%%%%%%%%%%%%%%%%%%%%%%%%%%%%%%%%%%%%%%%%%%%%%%%%%%%%%%%%%%%%%%%%%%%%%

\begin{myth}
Assume hypotheses (\ref{i}) - (\ref{l}). consider two  numbers $\alpha_1\neq 0$ and  $\alpha_2\neq 0$ and vectors
\begin{equation}
u^0 \in D(-\Delta)\cap H_{0}^{1}(\Omega),\, v_0 \in D(-\Delta),\;\mbox{and}\; u^1,\, v^1\in H_{0}^{1}(\Omega) \label{m}
\end{equation}
Then there exists a pair of functions $u,\, v $ in the class
\begin{equation}
\left\|\begin{array}{cc}
u,\,v\in L_{loc}^{\infty}(0,\,\infty;\,V)\\
u',\,v'\in L_{loc}^{\infty}(0,\,\infty;\,V) \label{n} \\
u^{''},\,v^{''}\in L_{loc}^{\infty}(0,\,\infty;\,L^{2}(\Omega))
\end{array}
\right.
\end{equation}
such that $u,\,v $ satisfy the system
\begin{eqnarray}
&&\displaystyle u^{''}(x,t)-mu(t)\Delta u(x,t)+\alpha_1\sum_{i=1}^{n}\frac{\partial v}{\partial x_i}=0,\; in\;L^{\infty}(0,\,\infty;\,L^{2}) \label{p}\\
&&\displaystyle v{''}(x,t)-\Delta v(x,t) -\alpha_2\sum_{i=1}^{n}\frac{\partial u}{\partial x_i}=0 ,\; in\;L^{\infty}(0,\,\infty;\,L^{2}) \label{q}
\end{eqnarray}
the boundary conditions
\begin{eqnarray}
&&\displaystyle \frac{\partial u}{\partial \nu} + h_1(.,u')=0\;\;in\, L^{\infty}_{loc}(0,\,\infty;\,L^1(\Gamma_1))\label{r}\\
&&\displaystyle \frac{\partial u}{\partial \nu} + h_2(.,v')+\sigma u=0\;\;in\, L^{\infty}_{loc}(0,\,\infty;\,L^1(\Gamma_1))\label{s}
\end{eqnarray}
and the initial conditions
\begin{equation}
u(0)=u^0,\;v(0)=v^0,\;u'(0)=u^1,\;v'(0)=v^1  \label{t}
\end{equation}
\label{teo1}
\end{myth}

In what follows, we introduce the notations and hypotheses to state the result on the decay of solutions. We will use the notations.
\begin{equation}
|u|\leq M||u||,\;\; ||u||_{L^2(\Gamma_1)}\leq N||u||,\;\;\forall u\in V \label{u}
\end{equation}

Consider the  function $ m(x)= x-x^0 , \; x\in \mathbb{R}^n$ ($x_0 $ a fixed vector de $\mathbb{R}^n$ ). Assume that there exist $x^0\in \mathbb{R}^n$ such that

\begin{equation}
\Gamma_0 = \{x\in\Gamma;\;m(x).\nu(x)\leq 0\}\;\;  \Gamma_1 = \{x\in\Gamma;\;m(x).\nu(x)> 0\}\label{w}
\end{equation}

Use the notations,
\begin{equation}
R(x^0)=max\{||m(x)||_{\mathbb{R}^n};\;x\in \overline{\Omega}\},\;\;0<\tau_0=min\{m(x).\nu(x);\; x\in\Gamma_1  \} \label{v}
\end{equation}

Assume the

\begin{equation}
h_1(x,r)=\Bigl[ m(x).\nu(x)\Bigr]p_1(s)\;\;h_2(x,r)=\Bigl[ m(x).\nu(x)\Bigr]p_2(s) \label{x}
\end{equation}

where $p_i\;\; (i=1,\,2)$  satisfy

\begin{equation}
\left\|\begin{array}{cc}
  p_i\in C^0(\mathbb{R}),\;p_i(0)=0\\
  \Bigl[p_i(s) -p_i(r)\Bigr](s-r)\geq b_i(s-r)^2,\;\;\forall s, r \in\mathbb{R} \label{y}\\
  |p_i(s)|\leq l_i|s|,\;\; s \in\mathbb{R}
\end{array}
\right.
\end{equation}
were $b_i,\;\;\mbox{and}\;\; L_i$ are positive constants.

We consider two real numbers $\alpha_1>0,\;\;\mbox{and}\;\;\alpha_2>0$. Introduce the following notations
\begin{equation}
A=2(n-1)\frac{M}{\mu^{\frac{1}{2}}_{0}}+2(n-1)M\frac{\alpha_1}{\alpha_2}+4\frac{R(x^0)}{\mu^{\frac{1}{2}}}+4R(x^0)\frac{\alpha_1}{\alpha_2}\label{a1}
\end{equation}

\begin{equation}
P_1=4(n-1)^2n\frac{M^2}{\mu_0}+16\mathbb{R}^2(x^0)\frac{n}{\mu_0}\label{a2}
\end{equation}

\begin{equation}
\left\|\begin{array}{cc}
P_2=4(n-1)^2||\displaystyle\sum_{i=1}^{n}\nu_i||_{L^{\infty}(\Gamma_1)}^{2}\frac{N^4}{\mu_0}+ \label{a3}\\
\frac{2N^2}{\tau_0\mu_0}||\displaystyle\sum_{i=1}^{n}\nu_i||^2 +16R^2(x^0)\frac{n}{\mu_0}
\end{array}
\right.
\end{equation}

\begin{equation}
S_1=4(n-1)^2\mu(0)R(x^0)L^{2}_{1}N^2+\mu(0)R^2(x^0)L^{2}_{2}+1 \label{a4}
\end{equation}

\begin{equation}
S_2=4(n-1)^2R(x^0)L^{2}_{2}N^2+2R^2(x^0)L^{2}_{2}+1 \label{b4}
\end{equation}

With respect to positive real numbers $\alpha_1,\;\;\mbox{and}\;\;\alpha_2$, we assume the following hypotheses:

\begin{equation}
\alpha_1\alpha_2\leq\frac{\mu_0}{64nM^2}\;\;  \mbox{and}\;\; P_1\alpha^2_1 + P_2\alpha^2_2\leq \frac{7}{8}  \label{a5}
\end{equation}

We consider a positive functions $\sigma(x)$ given by

\begin{equation}
\sigma(x)=\alpha_2\Bigl(\ \displaystyle\sum_{i=1}^{n}\nu_i(x)\Bigr)  \label{a6}
\end{equation}

We take three real numbers $\epsilon_1,\;\epsilon_2,\; \mbox{and}\;\; \eta$ satisfying

\begin{equation}
0< \epsilon_1 \leq \frac{1}{4A},\;\; 0<\epsilon_2 \leq min\{\frac{\mu_0b_1}{S_1},\; \frac{\bigl(\frac{\alpha_1}{\alpha_2}\bigr)b_2}{S_2}\} \label{a7}
\end{equation}

and

\begin{equation}
0 < \eta \leq min \{\epsilon_1,\;\epsilon_2 \} \label{a8}
\end{equation}

Introduction the energy

\begin{equation}
E(t)=\frac{1}{2}\Bigl\{|u'(t)|^2+\frac{\alpha_1}{\alpha_2}|v'(t)|^2+\mu(t)||u(t)||^2 +\frac{\alpha_1}{\alpha_2}||v(t)||^2 \Bigr \}, \; t\geq 0  \label{a9 }
\end{equation}

\begin{myth}
Assume hypotheses (\ref{w}), (\ref{x}), (\ref{y}), (\ref{a5}) and (\ref{a6}). Assume also that
\begin{equation}
\mu'(t)\leq 0,\;\;a.e.\, t\in(0,\infty) \label{a10}
\end{equation}
Then the pair of solutions u, v given in Theorem (2.1) satisfy
\begin{equation}
E(t)\leq 3E(0)e^{-\frac{2}{3}\eta t},\;\forall t>0  \label{a11}
\end{equation}
where $\eta > 0 $ was defined in (\ref{a8})
\label{teo2}
\end{myth}
%%%%%%%%%%%%%%%%%%%%%%%%%%%%%%%%%%%%%%%%%%%%%%%%%%%%%%%%%%%%%%%%%%%%%%%%%%%%%%%%%%%%%%%%%%%%%%%%%%%%%%%%%%%%%%%%%%%%%%%%%%%%%%%%%%%%%%%%%%%%%%%%%%%%%%%%%%%%%%
\section{Existence of Solutions}
\setcounter{equation}{0}
Before proving Theorem \ref{teo2}, we need of some previous results.

\begin{mylem}
Let $h(x,s)$ be a function satisfying the hypotheses (\ref{i}) and (\ref{j}) with $d_0 > 0$. Then there exists sequence $(h_l)$ of vectors of $C^0(\mathbb{R}; L^{\infty}(\Gamma_1))$ satisfying the following conditions:
\begin{eqnarray*}
&&(i)\; h_l(x,0)=0\;\;a.e.\;x\in\Gamma_1;\\[5pt]
&&(ii)\; \Bigl[h_l(x,s)-h_l(x,r)\Bigr](s-r)\geq d_0(s-r)^2,\;\;\forall s,r\in\mathbb{R}\;\; \mbox{and}\;\;a.e.\; \in\Gamma_1\\[5pt]
&&(iii)\; For \;any l\in \mathbb{N} \;there\; exists \;a \;function \;c_l\in L^{\infty}(\Gamma_1) \;satisfying\\[5pt]
&&\;\;\;\;\;\;\;\;|h_l(x,s)-h_l(x,s)|\leq c_l|s-r|, \;\forall s,r\in \mathbb{R}\;an \;a.e. \;in \; \Gamma_1\\[5pt]
&&(iv)\; (h_l)\;converges\;to \; h \;uniformly \; in \; bounded \; sets \; of \; \mathbb{R}, \; a.e.\; x\in \Gamma_1
\end{eqnarray*}
\label{L1}
\end{mylem}

%%%%%%%%%%%%%%%%%%%%%%%%%%%%%%%%%%%%%%%%%%%%%%%%%%%%%%%%%%%%%%%%%%%%%%%%%%%%%%%%%%%%%%%%%%%%%%%%%%%%%%%%%%%%%%%%%%%%%%%%%%%%%%%%%%%%%%%%%%%%%%%%%%%%%%%%%%%%%%%%%
\begin{mylem}
Let $ T > 0 $ be a real number. Consider the sequence $(w_l)$ of vectors of $L^2(0,T; H^{-\frac{1}{2}}(\Gamma_1))\cap L^1(0,T; L^1(\Gamma_1))$ and  $ w\in L^2(0,T; H^{-\frac{1}{2}}(\Gamma_1)) $
and $\chi \in L^1(0,T; L^1(\Gamma_1) $  such that
\begin{eqnarray*}
&&(i)\; w_l\rightarrow w \; weak \; in\; L^2(0,T; H^{-\frac{1}{2}}(\Gamma_1)) \\[5pt]
&&(ii)\;  w_l\rightarrow \chi  \; in \; L^1(0,T; L^1(\Gamma_1))
\end{eqnarray*}
then, $w=\chi$
\label{L2}
\end{mylem}
%%%%%%%%%%%%%%%%%%%%%%%%%%%%%%%%%%%%%%%%%%%%%%%%%%%%%%%%%%%%%%%%%%%%%%%%%%%%%%%%%%%%%%%%%%%%%%%%%%%%%%%%%%%%%%%%%%%%%%%%%%%%%%%%%%%%%%%%%%%%%%%%%%%
The proof of Lemma \ref{L1} can be found in Strauss \cite{Strauss} and Lemma \ref{L2}, in Louredo and Milla Miranda \cite{LM}.
%%%%%%%%%%%%%%%%%%%%%%%%%%%%%%%%%%%%%%%%%%%%%%%%%%%%%%%%%%%%%%%%%%%%%%%%%%%%%%%%%%%%%%%%%%%%%%%%%%%%%%%%%%%%%%%%%%%%%%%%%%%%%%%%%%%%%%%%%%%%%%%%%%%%
\subsection{Proof of Theorem 2.1}
Let $(h_{1l})$ and $(h_{2l})$ be two sequences in the conditions of Lemma \ref{L1} that approximate $h_1$ and $h_2$, respectively, consider two sequences $(u^1_l)$ and $(v^1_l)$ of vectors of $C_0^{\infty}(\Omega)$  such that
\begin{equation}
u^1_l\rightarrow u^1\;\; \mbox{and}\;\;v^1_l\rightarrow v^1\;\; in \;\; H_0^1(\Omega) \label{a12}
\end{equation}

Note that
\begin{equation}
\left\|\begin{array}{cc}
\frac{\partial u^0}{\partial \nu} + h_{1l}(.,u^1_l)=0 \; on \; \Gamma_1\;\forall l\label{a13}\\
\frac{\partial v^0}{\partial \nu} + h_{2l}(.,v^1_l) + \sigma u^0 =0 \; on \; \Gamma_1 \; \forall l
\end{array}
\right.
\end{equation}
Now, we fix  $l\in \mathbb{N}$  and construct a basis  $\{w_1^l,\,w_2^,\,w_3^l...\}\; of\; V\cap H^2(\Omega)$ such that $u^0,\; v^0,\; u^1_l,\; v_l^1$ belong to the subspace
 $\bigl[w_1^l,\,w_2^,\,w_3^l,\;w_4^l\bigl] $ spanned by the vectors $ w_1^l,\,w_2^,\,w_3^l,\;w_4^l $. Let $V_m = \bigl[w_1^l,\,w_2^l,\,w_3^l,\;...w_m^l\bigl]$ be the subspace of  $V\cap H^2(\Omega)$ spanned by $w_1^l,\,w_2^l,\,w_3^l,\;...w_m^l$.

{\bf Approximated Problem:} We find an approximate solution $u_{lm},\;v_{lm}$ of Problem (\ref{e}) - (\ref{g}) belonging to $V_m$, i. e.

 $$u_{lm}(t)=\displaystyle\sum_{i=1}^{n}g_{jlm}(t)w_j^l,\;\; v_{lm}(t)=\displaystyle\sum_{i=1}^{n}h_{jlm}(t)w_j^l$$

 and $u_{lm},\;v_{lm}$ is a solution of the system,

 \begin{equation}
 \left\|\begin{array}{cc}
 (u^{''}_{lm},\varphi) +\mu(( u_{lm},\varphi)) +\mu\int_{\Gamma_1} h_{1l}(.,u'_{lm})\varphi d\Gamma_1 + \\\alpha_1\displaystyle(\sum_{i=1}^{n}\frac{\partial v_{lm}}{\partial x_i},\varphi)=0, \forall
\varphi \in V_m\\ \label{a14}
(v^{''}_{lm},\psi) +(( v_{lm},\psi)) +\int_{\Gamma_1} h_{2l}(.,v'_{lm})\psi d\Gamma_1 + \int_{\Gamma_1} \sigma u_{lm}\psi d\Gamma_1 -\\\alpha_2\displaystyle(\sum_{i=1}^{n}\frac{\partial u_{lm}}{\partial x_i},\psi)=0, \forall\psi \in V_m\\
u_{lm}(0)=u^{0},\;v_{lm}(0)=v^{0},\; u'_{lm}(0)=u^{1}_{l},\;v'_{lm}(0)=v^{1}_{l}
\end{array}
\right.
 \end{equation}
 The above finite-dimensional system has a solution $u_{lm,\;v_{lm}}$ defined in $[0,\,t_{lm})$. The following estimates allow us to extend this solutions to the interval $[0, \infty]$.

\subsection{Estimates I}
 Consider $\varphi = u'_{lm}$ and $\psi =v'_{lm}$ in $(\ref{a13})_1$ and $(\ref{a13})_2$, respectively, we obtain

\begin{equation}
 \begin{array}{cc}
 \displaystyle\frac{d}{dt}\biggl[\frac{1}{2}|u'_{lm}|^2+\frac{\mu}{2}||u_{lm}||^2\biggr]+\mu\int_{\Gamma_1}h_1(.,u'_{lm})u'_{lm}d\Gamma_1 +\\
  \alpha_1\displaystyle(\sum_{i=1}^{n}\frac{\partial v_{lm}}{\partial x_i},u'_{lm})=\frac{\mu'}{2}||u_{lm}||^2 \label{a15}
\end{array}
\end{equation}

and

\begin{equation}
 \begin{array}{cc}
 \displaystyle\frac{d}{dt}\biggl[\frac{1}{2}|v'_{lm}|^2+\frac{1}{2}||v_{lm}||^2\biggr]+\displaystyle \int_{\Gamma_1}h_2(.,v'_{lm})v'_{lm}d\Gamma_1+\\ \displaystyle \int_{\Gamma_1} \sigma u_{lm}v'_{lm} d\Gamma_1- \alpha_2\displaystyle(\sum_{i=1}^{n}\frac{\partial u_{lm}}{\partial x_i},v'_{lm})=0 \label{a16}
\end{array}
\end{equation}

Introduce the notation
\begin{equation}
E_{lm}(t)=\frac{1}{2}|u'_{lm}(t)|^2+\frac{1}{2}|v'_{lm}(t)|^2 +\frac{\mu}{2}||u_{lm}(t)||^2 +||v_{lm}(t)||^2 \label{a17}
\end{equation}
we add the both sides of (\ref{a15}) and (\ref{a16}) and use hypothesis (\ref{k}) on $\mu$ and Lemma \ref{L1}, part(ii), applied to $h_{1l},\;h_{2l}$ we have

\begin{equation}
 \begin{array}{cc}
\frac{d}{dt}E_{lm}+\mu_0 d_1\displaystyle \int_{\Gamma_1}u'^{2}_{lm}d\Gamma_1 +d_2\displaystyle \int_{\Gamma_1}v'^{2}_{lm}d\Gamma_1\leq \frac{\mu'}{2}||u_{lm}||^2 -\\
\alpha_1\Bigl(\displaystyle \sum_{i=1}^{n}\frac{\partial v_{lm}}{\partial x_i},u'_{lm}\Bigr) -\displaystyle \int_{\Gamma_1}\sigma u_{lm}v'_{lm}d\Gamma_1  + \alpha_2\Bigl(\displaystyle
\sum_{i=1}^{n}\frac{\partial u_{lm}}{\partial x_i},v'_{lm}\Bigr)\label{a18}
\end{array}
\end{equation}
We find

$$\Biggl|\alpha_1\Biggl(\displaystyle \sum_{i=1}^{n}\frac{\partial v_{lm}}{\partial x_i}, u'_{lm}\Biggr)\Biggr|\leq \alpha_1^{2}n\Bigl(\frac{1}{2}||v_{lm}||^2\Bigr)+\frac{1}{2}|u'_{lm}|^2$$

Similarly

$$\Biggl|\alpha_2\Biggl(\displaystyle \sum_{i=1}^{n}\frac{\partial u_{lm}}{\partial x_i}, v'_{lm}\Biggr)\Biggr|\leq \alpha_2^{2}\frac{n}{\mu_0}\Bigl(\frac{\mu}{2}||u_{lm}||^2\Bigr)+\frac{1}{2}|v'_{lm}|^2$$

Also,

$$\Biggl|\displaystyle \int_{\Gamma_1}\sigma u_{lm}v'_{lm} d\Gamma_1 \Biggr|\leq \frac{N^2}{\mu_0 d_2}||\sigma||^2_{L^{\infty}(\Gamma_1)}\Bigl(\frac{\mu}{2}||u_{lm}||^2 \Bigr)+ \frac{d_2}{2}||v'||^2_{L^2(\Gamma_1)} $$

Taking into account the last three inequations in (\ref{a18}), derive

$$\frac{d}{dt}E_{lm}+\mu_0d_1\displaystyle \int_{\Gamma_1}u'^{2}_{lm}d\Gamma_1 +\frac{d_2}{2}\displaystyle \int_{\Gamma_1}v'^{2}_{lm}d\Gamma_1\leq \Bigl(|\mu'|+K\Bigr)E_{lm}$$

where $K$ is the constant

$$K=\alpha^2_1 + \alpha^2_2\frac{n}{\mu_0} +\frac{N^2}{\mu_0d_2}||\sigma||^2_{L^{\infty}(\Gamma_1)}$$

So, integrating the preceding inequality on $[0,\,t),\;\; t<t_{lm}$, we find
$$ E_{lm}+\mu_0d_1\displaystyle \int_{0}^{t}\int_{\Gamma_1}u'^{2}_{lm}d\Gamma_1 ds+\frac{d_2}{2}\displaystyle\int_{0}^{t} \int_{\Gamma_1}v'^{2}_{lm}d\Gamma_1ds\leq E_{lm}(0)+\displaystyle \int_{0}^{t}\Bigl(|\mu'|+K\Bigr)E_{lm}ds $$

Convergence (\ref{a12}) yield

$$E_{lm}(0)\leq \frac{1}{2}|u^1|^2+\frac{1}{2}|v^1|^2+\frac{\mu(0)}{2}||u^0||^2 +||v^0||^2 +1=L_0,\;\; \forall l\geq l_0 $$

Thus, the last two inequalities and Gronwall Lemma provide
\begin{equation}
 \begin{array}{cc}
\frac{1}{2}|u'_{lm}(t)|^2+\frac{1}{2}|v'_{lm}(t)|^2+\frac{\mu}{2}||u_{lm}(t)||^2+||v_{lm}(t)||^2+\mu_0d_1\displaystyle\int_{0}^{t} \int_{\Gamma_1}u'^{2}_{lm}d\Gamma_1ds+\label{a19}\\
\frac{d_2}{2}\displaystyle\int_{0}^{t} \int_{\Gamma_1}v'^{2}_{lm}d\Gamma_1ds\leq L_0exp\displaystyle\int_{0}^{t}(|\mu'|+K)dt=C(T),\;\;\forall t\in[0,\,T],\;\;\forall l\geq l_0
\end{array}
\end{equation}

where the constat $C(T) > 0$ is independent of $l\geq l_0$ and $m$. So

\begin{equation}
 \begin{array}{cc}
(u_{lm})\;\; and\;\; (v_{lm})\;\; are\;\; bounded\;\; in\;\; L^{\infty}_{loc}(0,\infty ; V) \label{a20}\\
(u'_{lm})\;\; and\;\; (v'_{lm})\;\; are \;\; bounded\;\; in\;\; L^{\infty}_{loc}(0,\infty ; L^2(\Omega))\\
(u'_{lm})\;\; and\;\; (v'_{lm})\;\; are\;\; bounded\;\; in\;\; L^{\infty}_{loc}(0,\infty ; L^2(\Gamma_1)))
\end{array}
\end{equation}

\subsection{Estimates II}
Differentiating with respect to $t$ the approximate equation $(\ref{a14})_1$  and making $\varphi= u^{''}_{lm}$ in the resulting expression, we obtain

\begin{eqnarray*}
&&\displaystyle\frac{d}{dt}\Bigl[\frac{1}{2}|u^{''}_{lm}|^2 +\frac{\mu}{2}||u^{'}_{lm}||^2 \Bigr] +\mu'((u_{lm},u^{''}_{lm})) +\mu'\displaystyle\int_{\Gamma_1}h_{1l}(.,u'_{lm})u^{''}_{lm}d\Gamma_1+\\
&&\displaystyle\mu\int_{\Gamma_1} h_{1l}(.,u'_{lm})(u^{''}_{lm})^2d\Gamma_1 + \alpha_1\Bigl(\displaystyle\sum_{i=1}^{n}\frac{\partial v'_{lm}}{\partial x_i} , u^{''}_{lm}\Bigr) = \frac{\mu'}{2}||u'_{lm}||^2
\end{eqnarray*}

Considering $\alpha= \frac{\mu'}{\mu}u^{''}_{lm} $ in $(\ref{a14})_1$ we get

\begin{eqnarray*}
&&\displaystyle \mu'((u_{lm},  u^{''}_{lm}))+ \mu'\displaystyle\int_{\Gamma_1}h_{1l}(.,u'_{lm})u^{''}_{lm}d\Gamma_1 = -\frac{\mu'}{\mu}|u'_{lm}|^2-\frac{\mu'}{\mu}\Bigl(\alpha_1\displaystyle
\sum_{i=1}^{n}\frac{\partial v_{lm}}{\partial x_i}, u^{''}_{lm}\Bigl)\\
\end{eqnarray*}

Combining the last two equations, we final
\begin{eqnarray*}
&&\displaystyle\frac{d}{dt}\Bigl[\frac{1}{2}|u^{''}_{lm}|^2 +\frac{\mu}{2}||u^{'}_{lm}||^2 \Bigr]  +\mu\displaystyle\int_{\Gamma_1}h'_{1l}(.,u'_{lm})(u^{''}_{lm})^2d\Gamma_1= \frac{\mu'}{2}||u'_{lm}||^2 +\\
&&\displaystyle\frac{\mu'}{\mu}|u^{''}_{lm}|^2 + \alpha_1\frac{\mu'}{\mu}\Bigl(\displaystyle \sum_{i=1}^{n}\frac{\partial v_{lm}}{\partial x_i}, u^{''}_{lm} \Bigr) -\alpha_1 \Bigl(\displaystyle \sum_{i=1}^{n}\frac{\partial v'_{lm}}{\partial x_i},u^{''}_{lm} \Bigr)
\end{eqnarray*}

In similar way, approximate equation $(\ref{a14})_2$ provide,

 \begin{eqnarray*}
&&\displaystyle\frac{d}{dt}\Bigl[\frac{1}{2}|v^{''}_{lm}|^2 +\frac{1}{2}||v^{'}_{lm}||^2 \Bigr]  +\displaystyle\int_{\Gamma_1}h'_{2l}(.,v'_{lm})(v^{''}_{lm})^2d\Gamma_1= \\
&&\alpha_2\Bigl(\displaystyle \sum_{i=1}^{n}\frac{\partial u'_{lm}}{\partial x_i}, v^{''}_{lm} \Bigr) -\displaystyle \int{\Gamma_1}\sigma(x)u'_{lm}v^{''}d\Gamma
\end{eqnarray*}

Introduce the notation
 $$E^*_{lm}(t)=\frac{1}{2}|u^{''}_{lm}|^2+ \frac{1}{2}|v^{''}_{lm}|^2+\frac{\mu}{2}||u^{'}_{lm}||^2+ \frac{1}{2}|v^{'}_{lm}|^2,\;\; t\geq0$$

 Adding the both sides of the las two equations and using hypothesis (\ref{k}) on $\mu$ an Lemma \ref{L1}, part $(ii)$, applied to $h_{1l}$,  $h_{2l}$, we get,

\begin{eqnarray}
&&\displaystyle \frac{d}{dt}E^{*}_{lm}  +\mu d_1\displaystyle\int_{\Gamma_1}(u^{''}_{lm})2d\Gamma_1+ d_2\displaystyle\int_{\Gamma_1}(v^{''}_{lm})^2d\Gamma_1\leq \frac{\mu'}{2}||u'_{lm}||^2+ \nonumber\\
&&\displaystyle \alpha_1\frac{\mu'}{\mu}\Bigl(\displaystyle \sum_{i=1}^{n}\frac{\partial v_{lm}}{\partial x_i}, u^{''}_{lm} \Bigr) -\alpha_1 \Bigl(\displaystyle \sum_{i=1}^{n}\frac{\partial v'_{lm}}{\partial x_i},u^{''}_{lm} \Bigr) +\alpha_2 \Bigl(\displaystyle \sum_{i=1}^{n}\frac{\partial u'_{lm}}{\partial x_i},v^{''}_{lm} \Bigr)-\label{a21}\\
&&\displaystyle \int_{\Gamma_1}\sigma u'_{lm}v^{''}_{lm} d\Gamma_1 \nonumber
\end{eqnarray}

we have

$$\begin{array}{l}
{\rm \bullet} \;\;\;{\displaystyle \Bigl|\displaystyle \alpha_1\frac{\mu'}{\mu}\Bigl(\displaystyle \sum_{i=1}^{n}\frac{\partial v_{lm}}{\partial x_i}, u^{''}_{lm} \Bigr)  \Bigr|\leq \frac{(\alpha_1)^2}{\mu_{0}^{2}}n|\mu'|(\frac{1}{2}|v_{lm}|^2)+|\mu'|(\frac{1}{2})|u^{''}_{lm}|^2 }\\[5pt]
{\rm \bullet}\;\;\; {\displaystyle  \displaystyle \Bigl|\displaystyle \alpha_1\Bigl(\displaystyle \sum_{i=1}^{n}\frac{\partial v'_{lm}}{\partial x_i}, u^{''}_{lm} \Bigr)  \Bigr|\leq
(\alpha_1)^2n(\frac{1}{2}||v'_{lm}||^2)+\frac{1}{2}|u^{''}_{lm}|^2}\\[5pt]
{\rm \bullet}\;\;\; {\displaystyle \Bigl|\alpha_2 \Bigl(\displaystyle\sum_{i=1}^{n}\frac{\partial u'_{lm}}{\partial x_i} ,v^{''}_{lm}\Bigr) \Bigr|\leq \frac{\alpha^2_2}{\mu_0}n(\frac{\mu}{2})||u'_{lm} ||^2+\frac{1}{2}|v^{''}_{lm}|^2 }\\[5pt]
{\rm \bullet}\;\;\; {\displaystyle \Bigl|\displaystyle \int_{\Gamma_1}\sigma u'_{lm}v^{''}d\Gamma_1 \Bigr| \leq \frac{N^2}{d_2\mu_0}||\sigma||^2_{L^{\infty}(\Gamma_1)}(\frac{\mu}{2}||u^{2}_{lm}||^2) + \frac{1}{2}||v^{''}_{lm}||^{2}_{L^2(\Gamma_1)}}
\end{array}$$
from the last four inequalities in (\ref{a21}) and using the boundedness (\ref{a19})  for $||v_{lm}||^2$ we have
\begin{equation}
\frac{d}{dt}E^*_{lm} + \mu_0d_1\displaystyle\int_{\Gamma_1}(u_{lm}^{''})^2d\Gamma_1 + \frac{d_2}{2}\displaystyle\int_{\Gamma_1}(v_{lm}^{''})^2d\Gamma_1 \leq K_1(T)|\mu'|+\Bigl[K_2|\mu'| +K_3\Bigr]E^*_{lm}\label{a22}
\end{equation}
where
\begin{equation}
K_1(T)=\frac{\alpha^2_1}{\mu_0}n C(T),\;\;K_2=\frac{1}{\mu_0}+1,\;\; K_3 = \alpha^2_1n+ \frac{\alpha^2_2}{\mu_0}n+\frac{N^2}{d_2\mu_0}||\sigma||_{L^{\infty}(\Gamma_1)} +2 \label{a23}
\end{equation}

Integrate both sides of (\ref{a22}) on $[0,\,t],\;\;0<t\leq T$, we find

\begin{equation}
E^*_{lm}(t)\leq E^*_{lm}(0) + K_1(T)\displaystyle \int_{0}^{T}|\mu'|dt + \displaystyle \int_{0}^{T}[K_2|\mu'|+k_3]E^*_{lm}ds \label{a24}
\end{equation}

We will obtain a second estimate if we bound $E^*_{lm}(0)$. This is the key point of the proof of theorem \ref{teo1}. The boundedness will follow by the choice of the special basis of $V\cap H^2(\Omega).$

In fact, if we make $t=0$ in the approximate equations $(\ref{a14})_1$ and $(\ref{a14})_2$ and consider $\varphi = u''_{lm}$ and $\psi = u''_{lm}$ we have
\begin{equation}
|u^{''}_{lm}(0)|^2+\mu(0)((u^0, u^{''}_{lm}(0)))+\mu(0)\displaystyle\int_{\Gamma_1}h_{1l}(.,u^1_l)u^{''}_{lm}d\Gamma_1+\alpha_1\Bigl(\displaystyle\sum_{i=1}^{n}\frac{\partial v^0}{\partial x_i}, u^{''}_{lm}(0)\Bigr) =0 \label{a25}
\end{equation}
and
\begin{eqnarray}
&&|v^{''}_{lm}(0)|^2+\mu(0)((v^0, v^{''}_{lm}(0)))+\displaystyle\int_{\Gamma_1}h_{2l}(.,v^1_l)v^{''}_{lm}d\Gamma_1+\displaystyle \int_{\Gamma_1}\sigma u^0u^{''}_{lm}d\Gamma_1+\nonumber\\
&&\alpha_1\Bigl(\displaystyle\sum_{i=1}^{n}\frac{\partial v^0}{\partial x_i}, u^{''}_{lm}(0)\Bigr) =0. \label{a26}
\end{eqnarray}

The Gauss Theorem and the equalities (\ref{a13}) provide
\begin{eqnarray*}
 &&\displaystyle\mu(0)((u^0,\, u^{''}_{lm}(0)))+\mu(0)\displaystyle \int_{\Gamma_1}h_{1l}(.,u^1_l)u^{''}_{lm}(0)d\Gamma_1=\\
 &&\displaystyle\mu(0)\Bigl[(-\Delta u^0, u^{''}_{lm}(0)) +\int_{\Gamma_1} (\frac{\partial u^0}{\partial \mu} + h{1l}(., u^1_l))u^{''}_{lm}(0)d \Gamma_1\Bigr]=\\
 &&\displaystyle\mu(0)(-\Delta u^0, u^{''}_{lm}(0)).
 \end{eqnarray*}

Taking into account in the last two equations  (\ref{a25}) and (\ref{a26}), we get
 $$|u^{''}_{lm}(0)|^2+\mu(0)(-\Delta u^0, u^{''}_{lm}(0)) + \alpha_1\Bigl(\displaystyle \sum_{i=1}^{n}\frac{\partial v^0}{\partial x_i}, u^{''}_{lm}(0) \Bigr) =0 $$
 and
$$|v^{''}_{lm}(0)|^2+\mu(0)(-\Delta v^0, v^{''}_{lm}(0)) - \alpha_2\Bigl(\displaystyle \sum_{i=1}^{n}\frac{\partial u^0}{\partial x_i}, v^{''}_{lm}(0) \Bigr) =0. $$
So
\begin{equation}
|u^{''}_{lm}(0)|^2 \leq \mu(0)|\Delta u^0| + |\alpha_1 |n^{\frac{1}{2}}||v^0|| = a_1 \label{a27}
 \end{equation}
 and
 \begin{equation}
|v^{''}_{lm}(0)|^2 \leq \mu(0)|\Delta v^0| + |\alpha_2 |n^{\frac{1}{2}}||u^0|| = a_2. \label{a28}
 \end{equation}

Therefore, the last two boundedness and convergence (\ref{a12}) provide
\begin{equation}
E^*_{lm}(0)\leq \frac{a_1^2}{2}+ \frac{a_2^2}{2}+\frac{\mu(0)}{2}||u^1||^2+ \frac{1}{2}||v^1 ||^2 +1 = a_3,\; \forall l\geq l_0. \label{a29}
\end{equation}

The inequalities (\ref{a27}) and (\ref{a29}) and Gronwall Lemma yields
\begin{equation}
 \begin{array}{cc}
\frac{1}{2}|u^{''}_{lm}(t)|^2 +\frac{1}{2}|v^{''}_{lm}(t)|^2 +\frac{\mu(t)}{2}||u'_{lm}(t)||^2 +\frac{1}{2}||v'_{lm}(t)||^2 +\mu_0d_1\displaystyle \int_{\Gamma_1}(u^{''}_{lm})^2d\Gamma_1 \\ +\frac{d_2}{2}\displaystyle \int_{\Gamma_1}(v^{''}_{lm})^2d\Gamma_1 \leq \Bigl[a_3+K_1(T)\displaystyle \int_0^t|\mu'|dt\Bigr]exp\displaystyle\int_0^t\Bigl[K_2|\mu'|+k_3\Bigr]dt=C_1(T)\\
\forall t\in [0, T],\;\; l\geq l_0 \label{a30}
\end{array}
\end{equation}
where $C_1(T) > 0$ is a constant independent of $l\geq l_0$ and $m. $

Thus
\begin{equation}
 \left\|\begin{array}{cc}
 (u'_{lm})\;\mbox{and}\; (v'_{lm})\;\mbox{are bounded in } L^{\infty}_{loc}(0,\infty ; V)\\
 (u^{''}_{lm})\;\mbox{and}\; (v^{''}_{lm})\;\mbox{are bounded in } L^{\infty}_{loc}(0,\infty ; L^2(\Omega))\label{a31} \\
 (u^{''}_{lm})\;\mbox{and}\; (v^{''}_{lm})\;\mbox{are bounded in } L^{2}_{loc}(0,\infty ; L^2(\Gamma_1))
\end{array}
\right.
 \end{equation}

The boundedness (\ref{a20}) and (\ref{a31}) provide two subsequences of $(u_{lm})$ and $(v_{lm})$, still denoted by $u_{lm}$ and $v_{lm}$, and two function $u_l$ and $v_l$ such that
\begin{equation}
 \left\|\begin{array}{cc}
 u_{lm}\rightarrow u_l\,\, \mbox{and}\,\,v_{lm}\rightarrow v_l \,\,\mbox{werk star in}\,\, L^{\infty}_{loc}(0,\infty ; V)\\
 u'_{lm}\rightarrow u'_l\,\, \mbox{and}\,\,v'_{lm}\rightarrow v'_l \,\,\mbox{werk star in}\,\, L^{\infty}_{loc}(0,\infty ; V)\\
 u^{''}_{lm}\rightarrow u^{''}_l\,\, \mbox{and}\,\,v^{''}_{lm}\rightarrow v^{''}_l \,\,\mbox{werk star in}\,\, L^{\infty}_{loc}(0,\infty ; L^2(\Omega))\label{a32}\\
 u'_{lm}\rightarrow u'_l\,\, \mbox{and}\,\,v'_{lm}\rightarrow v'_l \,\,\mbox{werk in}\,\, L^{2}_{loc}(0,\infty ; L^2(\Gamma_1))\\
 u^{''}_{lm}\rightarrow u^{''}_l\,\, \mbox{and}\,\,v^{''}_{lm}\rightarrow v^{''}_l \,\,\mbox{werk in}\,\, L^{2}_{loc}(0,\infty ; L^2(\Gamma_1))\\
\end{array}
\right..
 \end{equation}

\subsection{Passage to the Limit in m}

  We analyze the nonlinear terms on the boundary $\Gamma_1$. Let $T > 0$ be a real number. By convergence $(\ref{a32})_2$ and $(\ref{a32})_5$, the compact embedding of $H^{\frac{1}{2}}(\Gamma_1)$ in $L^2(\Gamma_1)$ and the Aubin-Lions Theorem [9 ], give us

  $$ u^{'}_{lm} \rightarrow u'_l\,\, \mbox{in}\,\, L^{2}_{loc}(0,\infty ; L^2(\Gamma_1))$$

  Lemma \ref{L1}, part (iv), provide

  $$\displaystyle\int_{\Gamma_1}|h_{1l}(., u'_{lm}) - h_{1l}(., u'_{l})|^2d\Gamma_1\leq ||c_l||_{L^{\infty}(\Gamma_1)}\displaystyle\int_{\Gamma_1}|u'_{lm}-u'_{l}|^2d\Gamma_1$$

These two results yield

$$h_{1l}(.,u'_{lm}) \rightarrow  h_{1l}(.,u'_{l}) \,\,\mbox{in}\,\, L^{2}_{loc}(0,\infty ; L^2(\Gamma_1))$$

Then by a diagonal process, we obtain

\begin{equation}
 h_{1l}(.,u'_{lm}) \rightarrow  h_{1l}(.,u'_{l}) \,\,\mbox{in}\,\, L^{2}_{loc}(0,\infty ; L^2(\Gamma_1))\label{a33}
\end{equation}

In a similar way, we find

\begin{equation}
 h_{2l}(.,v'_{lm}) \rightarrow  h_{2l}(.,v'_{l}) \,\,\mbox{in}\,\, L^{2}_{loc}(0,\infty ; L^2(\Gamma_1))\label{a34}
 \end{equation}

We take the limit in $m$ of system (\ref{a14}). Then by convergence (\ref{a32}),  (\ref{a33}),  (\ref{a34}) and noting that $V_m$ is dense in $V$, we obtain

\begin{equation}
 \begin{array}{cc}
\displaystyle\int_{0}^{\infty}(u^{''}_l,\varphi)\theta dt + \displaystyle\int_{0}^{\infty}\mu((u_l,\varphi))\theta dt +\displaystyle\int_{0}^{\infty}\int_{\Gamma_1}\mu h_{1l}(.,u'_l)\varphi\theta d\Gamma_1dt +\\
\displaystyle\int_{0}^{\infty}\alpha_1\Bigl(\displaystyle \sum_{i=1}^{n}\frac{\partial v_l}{\partial x_i}, \varphi\Bigr)dt = 0\;;\ \forall \varphi \in V, \forall \theta \in D(0, \infty) \label{a35}
\end{array}
\end{equation}
and

\begin{equation}
 \begin{array}{cc}
\displaystyle\int_{0}^{\infty}(v^{''}_l,\psi)\theta dt + \displaystyle\int_{0}^{\infty}((v_l,\psi))\theta dt +\displaystyle\int_{0}^{\infty}\int_{\Gamma_1} h_{2l}(.,v'_l)\psi\theta d\Gamma_1dt +\\
\displaystyle\int_{0}^{\infty}\int_{\Gamma_1}\sigma u_l \psi \theta - \displaystyle\int_{0}^{\infty}\alpha_2\Bigl(\displaystyle \sum_{i=1}^{n}\frac{\partial u_l}{\partial x_i}, \psi\Bigr)dt = 0\;;\ \forall \psi \in V, \forall \theta \in D(0, \infty) \label{a36}
\end{array}
\end{equation}

Considering $\varphi$, $\psi$ in $D(\Omega)$ in the preceding equations and noting the regularity of $u_l$, $v_l$ given in (ref{a32}), we get

\begin{equation}
 \begin{array}{cc}
u^{''}_l-\mu \Delta u_l +\alpha_1 \displaystyle \sum_{i=1}^{n}\frac{\partial v_l}{\partial x_i}=0 \;\; in \;\;L^{\infty}_{loc}(0, \infty; L^2(\Omega))\\
v^{''}_l-\Delta v_l - \alpha_2\displaystyle\sum_{i=1}^{n}\frac{\partial u_l}{\partial x_i}=0 \;\; in\;\; L^{\infty}_{loc}(0, \infty; L^2(\Omega)) \label{a37}
 \end{array}
\end{equation}

This implies that $\Delta u_l$, $\Delta v_l$ belongs to $L^{\infty}_{loc}(0, \infty ; L^2(\Omega))$ and $u_l$ and $v_l$ belong to $L^{\infty}_{lon}(0, \infty ; V)$, we find
$\displaystyle\frac{\partial u_l}{\partial x_i} $, $\displaystyle\frac{\partial v_l}{\partial x_i}$ in $L^2_{loc}(0, \infty; H^{-\frac{1}{2}}(\Gamma_1))$,
%%%%%%%%%%%%%%%%%%%%%%%%%%%%%%%%%%%%%%%%%%%%%%%%%%%%%%%%%%%%%%%%%%%%%%%%%%%%%%%%%%%%%%%%%%%%%%%%%%%%%%%%%%%%%%%%%
see \cite{Milla}.
%%%%%%%%%%%%%%%%%%%%%%%%%%%%%%%%%%%%%%%%%%%%%%%%%%%%%%%%%%%%%%%%%%%%%%%%%%%%%%%%%%%%%%%%%%%%%%%%%%%%%%%%%%%%%%%%%%

Multiplying both sides of equation (\ref{a37}) by $\varphi\theta$ and  $\psi\theta$ with $\varphi,\;\;\psi$ in $V$ and $\theta \in D(0, \infty)$, using the green formulae and preceding regularity, we obtain

\begin{equation}
\begin{array}{cc}
\displaystyle\int_{0}^{\infty}(u^{''}_l, \varphi)\theta dt +\displaystyle\int_{0}^{\infty}\mu ((u_l, \varphi))\theta dt -\displaystyle\int_{0}^{\infty}<\mu\frac{\partial u_l}{\partial \nu}, \varphi>\theta dt + \\ \displaystyle \int_{0}^{\infty}\alpha_1 \Bigl(\sum_{i=1}^{n}\frac{\partial v_l}{\partial x_i},\varphi\Bigr)\theta dt \label{a38}\\
\end{array}
\end{equation}

and

\begin{equation}
\begin{array}{cc}
\displaystyle\int_{0}^{\infty}(v^{''}_l, \psi)\theta dt +\displaystyle\int_{0}^{\infty}((v_l, \psi))\theta dt -\displaystyle\int_{0}^{\infty}<\frac{\partial v_l}{\partial \nu}, \psi>\theta dt + \\ \displaystyle \int_{0}^{\infty}\alpha_2 \Bigl(\sum_{i=1}^{n}\frac{\partial u_l}{\partial x_i},\psi\Bigr)\theta dt \label{a39}\\
\end{array}
\end{equation}

where $<.\, ,\, .> $ is the duality pairing between $H^{-\frac{1}{2}}(\Gamma_1)$ and  $H^{\frac{1}{2}}(\Gamma_1)$ . Comparing equations (\ref{a38}) and (\ref{a35}) with (\ref{a39}) and (\ref{a36}),  using the regularity of $h_{1l}(.,u'_l)$ and $h_{2l}(.,v'_l)$ given in (\ref{a33}) and (\ref{a34}), respectively, we have

\begin{equation}
\begin{array}{cc}
 \displaystyle \frac{\partial u_l}{\partial \nu} + h_{1l}(., u'_l) =0 \;\; in\;\; L^{\infty}_{loc}(0, \infty , L^2(\Gamma_1))\label{a40} \\
 \displaystyle \frac{\partial v_l}{\partial \nu} + h_{2l}(., v'_l) + \sigma u_l =0 \;\; in\;\; L^{\infty}_{loc}(0, \infty , L^2(\Gamma_1))
\end{array}
\end{equation}

\subsection{Passage to the Limit in l}

As the boundedness (\ref{a19}) and (\ref{a30}) are independent of $l\geq l_0$ and $m$, we obtain analogous convergence to (\ref{a18}), i. e., there are functions $u$ and $v$ such that

\begin{equation}
 \left\|\begin{array}{cc}
 u_{l}\rightarrow u \,\, \mbox{and}\,\,v_{l}\rightarrow v \,\,\mbox{werk star in}\,\, L^{\infty}_{loc}(0,\infty ; V)\\
 u'_{l}\rightarrow u'\,\, \mbox{and}\,\,v'_{l}\rightarrow v'_l \,\,\mbox{werk star in}\,\, L^{\infty}_{loc}(0,\infty ; V)\\
 u^{''}_{l}\rightarrow u^{''}\,\, \mbox{and}\,\,v^{''}_{l}\rightarrow v^{''} \,\,\mbox{werk star in}\,\, L^{\infty}_{loc}(0,\infty ; L^2(\Omega))\label{a41}\\
 u'_{l}\rightarrow u'\,\, \mbox{and}\,\,v'_{l}\rightarrow v' \,\,\mbox{werk in}\,\, L^{2}_{loc}(0,\infty ; L^2(\Gamma_1))\\
 u^{''}_{l}\rightarrow u^{''}\,\, \mbox{and}\,\,v^{''}_{l}\rightarrow v^{''} \,\,\mbox{werk in}\,\, L^{2}_{loc}(0,\infty ; L^2(\Gamma_1))\\
\end{array}
\right.
 \end{equation}

These convergence allow us to pass to the limit in (\ref{a35}) and (\ref{a36}). So for $ \varphi$, $\psi$ in $D(\Omega)$, we obtain

\begin{equation}
 \left\|\begin{array}{cc}
 u^{''} - \mu\Delta u + \alpha_1\displaystyle\sum_{i=1}^{n}\frac{\partial v}{\partial x_i} =0 \;\; L^{\infty}_{loc}(0, \infty , L^2(\Gamma_1))\\
 v^{''} - \mu\Delta v + \alpha_2\displaystyle\sum_{i=1}^{n}\frac{\partial u}{\partial x_i} =0 \;\; L^{\infty}_{loc}(0, \infty , L^2(\Gamma_1)) \label{a42}
\end{array}
\right.
\end{equation}
In what follow, we analyze the equation (\ref{a40}). Let $T > 0$ be a fixed real number. The convergence $(\ref{a41})_2$ yields
$$u'_l\rightarrow u'\;\; wark \;\;in \;\; L^2(\Gamma_1)$$
This, the compact immersion of $H^{\frac{1}{2}}(\Gamma_1)$ in $L^2(\Gamma_1)$ and the Aubin - Lions Theorem, give us
$$u'_l \rightarrow u'\;\; in\;\; L^{2}(0, T: L^2(\Gamma_1)) $$
which implies
$$u'_l(x,t) \rightarrow u'(x,t)\;\;a.e.\;\; x\in \Gamma_1,\;\;l\in (0,T). \label{a43} $$

Analogously,
 $$v'_l(x,t) \rightarrow v'(x,t)\;\;a.e.\;\; x\in \Gamma_1,\;\;l\in (0,T). \label{a44} $$

Fix $(x,t)\in \Gamma_1 \times (0, T)$. The last convergence implies that the set $\{u'_l(x,t),\;\;v'_l(x,t);\;\; l\geq l_0\}$ is a bounded set of $\mathbb{R}$. This and Lemma \ref{L1}, part (iv), on the uniformly convergence of $h_{1l}(x,u'(x,t))$ and  $h_{2l}(x,u'(x,t))$, provide
\begin{equation}
\begin{array}{cc}
h_{1l}(x,u'_l(x,t))\rightarrow h_1(x,u'(x,t))\;\; in\;\;a.e.\;\;  \Gamma_1\times (0, T)\\
h_{2l}(x,v'_l(x,t))\rightarrow h_2(x,v'(x,t))\;\; in\;\; \Gamma_1\times (0, T) \label{a45}
\end{array}
\end{equation}

We take the scalar product of $L^2(\Omega)$ on both sides of equation $(\ref{a37})_1$ with $u'_l$ and integrate on $[0, T]$  to obtain

\begin{eqnarray*}
&&\displaystyle\int_0^{T}\int_{\Gamma_1}\mu h_{1l}(.,u'_l)d\Gamma_1 dt =-\frac{1}{2}|u'_l(T)|^2 +\frac{1}{2}|u'_l(0)|^2-\frac{\mu(T)}{2}||u_l(T)||^2+\frac{\mu(0)}{2}||u^0||^2 - \\
&&\displaystyle\int_0^{T}\alpha_1\Bigl(\sum_{i=1}^{n}\frac{\partial v_l}{\partial x_i}, u'_l\Bigr) + \displaystyle\int_0^{T}\frac{\mu'}{2}||u_l||^2 dt
\end{eqnarray*}

By estimate (\ref{a19}) we find that each term of the second member of the preceding expression can be bound by a constant $C_3(T)$. Thus

 \begin{eqnarray}
 \displaystyle\int_0^{T}\int_{\Gamma_1}h_{1l}(.,u'_l)u'_ld\Gamma_1dt\leq \frac{C_4(T)}{\mu_0} \label{a46}
\end{eqnarray}

In a similar way, we get from $(\ref{a37})_2$ that

\begin{eqnarray}
 \displaystyle\int_0^{T}\int_{\Gamma_1}h_{2l}(.,v'_l)v'_ld\Gamma_1dt\leq C_5(T) \label{a47}
\end{eqnarray}

The constants $C_4(T)$ and $C_5(T)$ are independent of $l\geq l_0$. The results (\ref{a45}) - (\ref{a47}) allow us to apply the Strauss Theorem \cite{Strauss} to obtain

\begin{equation}
\begin{array}{cc}
h_{1l}(.,u'_l)\rightarrow h_1(.,u')\;\; in\;\; L^1(\Gamma_1\times (0, T))\\
h_{2l}(.,v'_l)\rightarrow h_2(.,v')\;\; in\;\; L^1(\Gamma_1\times (0, T)) \label{a48}
\end{array}
\end{equation}

On the other hand, by convergence (\ref{a41}) and equation $(\ref{a37})_1$ we deduce that
\begin{eqnarray*}
&& u_j\rightarrow u \;\; weak \;\;in\;\; L^2(0, T; V) \;\; and \;\; \Delta u_l\rightarrow \Delta u \; weak \;\;in\;\; L^2(0, T; L^2(\Omega)).
\end{eqnarray*}

Therefore
\begin{eqnarray*}
&& \displaystyle \frac{\partial u_l}{\partial\nu} \rightarrow \displaystyle \frac{\partial u}{\partial\nu} \;\;weak\;\; in L^2(o, T; H^{-\frac{1}{2}}(\Gamma_1))
\end{eqnarray*}

See \cite{Milla}. By equation $(\ref{a40})_1$ and convergence $(\ref{a48})_1$, we get

\begin{eqnarray*}
&& \displaystyle \frac{\partial u_l}{\partial\nu} = - h_{1l}(.,u'_l)\rightarrow h_1(., u')\;\; in\;\; L^1(0, T; L^1(\Gamma_1))
\end{eqnarray*}

The last convergence and Lemma \ref{L2} provide

\begin{eqnarray*}
&& \displaystyle \frac{\partial u_l}{\partial\nu} + h_1(.,u'_l)=0\;\; in\;\; L^1(0, T; L^1(\Gamma_1))
\end{eqnarray*}

Then by a diagonal process, we obtain

 \begin{eqnarray}
&& \displaystyle \frac{\partial u_l}{\partial\nu} + h_1(.,u'_l)=0\;\; in\;\; L^1_{loc}(0, \infty; L^1(\Gamma_1))\label{a49}
\end{eqnarray}

In a similar way, we deduce

\begin{eqnarray}
&& \displaystyle \frac{\partial v_l}{\partial\nu} + h_2(.,v'_l)+ \sigma u=0\;\; in\;\; L^1_{loc}(0, \infty; L^1(\Gamma_1))\label{a50}
\end{eqnarray}

The initial condition (\ref{t}) is obtained from $(\ref{a14})_3$ and the estimates (\ref{a32}), (\ref{a41}).

With the above part, the equations (\ref{a42}), (\ref{a48}), (\ref{a50}) and the estimates (\ref{a41}) we get the proof of theorem \ref{teo1}.

\section{Decay of Solutions}
\setcounter{equation}{0}
 Before proving the theorem \ref{teo2}, we introduce some previous results.

 \begin{mypro}
Let $h:\mathbb{R} \rightarrow \mathbb{R} $ be a Lipschitz continuous function. If $u\in H^{\frac{1}{2}}(\Gamma_1)$, then $h(u)\in H^{\frac{1}{2}}(\Gamma_1)$ and the map $h:H^{\frac{1}{2}}(\Gamma_1)\rightarrow L^{\frac{1}{2}}(\Gamma_1) $ is continuous.
\label{P4}
\end{mypro}

Let $(p_{il})$ be the sequence of Lipschitz continuous functions given in Lema \ref{L1} that approximate $p_i\;\;(i=1,\,2)$. Note that $u'_l\in L^{\infty}_{loc}(0,\infty;H^{\frac{1}{2}}(\Gamma_1))$ $(see  (\ref{a32})_2)$  and $\Gamma$ is of class $C^2$. The proposition \ref{P4} implies that
$$(m\cdot\nu)p_{1l}(u'_l)\in L^{\infty}_{loc}(0,\infty;H^{\frac{1}{2}}(\Gamma_1)). $$

This and $(\ref{a37})_1$ provide
$$\frac{\partial u_l}{\partial \nu} =-(m\cdot \nu)p_{1l}(u'_l)=g_l \in L^{\infty}_{loc}(0,\infty;H^{\frac{1}{2}}(\Gamma_1)) $$
Also,the  equation $(\ref{a37})_1$ implies that $\Delta u_l \in L^{\infty}_{loc}(0,\infty;L^2(\Omega)).$

Thus $u_l(t)$ is the solution of the following elliptic problem
\begin{eqnarray*}
-\Delta u_l(t)=f_t\;\; in\;\;\Omega\;\; (f_l(t)\ in L^2(\Omega)) \\
u_l(t) = 0\;\; on \;\;\Gamma_0 \;\;\;\;\;\;\;\;\;\;\;\;\;\;\;\;\;\;\;\;\;\;\;\;\;\;\;\;\;\;\;\\
\frac{\partial u_l(t)}{\partial \nu}= g_l\;\; on \;\; \Gamma_1 \;\; (g_l(t))\in H^{\frac{1}{2}}(\Gamma_1)
\end{eqnarray*}

By regularity of elliptic problems we have
\begin{eqnarray}
&&\displaystyle u_l \in L^{\infty}_{loc}(0, \infty; V\cap H^2(\Omega)) \label{D1}
\end{eqnarray}
(see \cite{Milla}). Similarly,
\begin{eqnarray}
&&\displaystyle v_l \in L^{\infty}_{loc}(0, \infty; V\cap H^2(\Omega)).\label{D2}
\end{eqnarray}

The regularity (\ref{D1}) allows us to obtain the following identities
 \begin{equation}
\left\|\begin{array}{cc}
(\Delta u_l, m\cdot \nabla u_l)=(n-2)|| u_l||^2 - \displaystyle \int_{\Gamma}(m\cdot \nu)|\nabla u_l|^2d\Gamma + 2\displaystyle \int_{\Gamma}\frac{\partial u_l}{\partial \nu} (m\cdot \nabla u_l)d\Gamma\\
2(u'_l,m\cdot \nabla u'_l)=-n|u'_l|^2 +\displaystyle \int_{\Gamma}(m\cdot \nu)(u')^{2}d\Gamma \label{D3}
\end{array}
\right.
\end{equation}
(see \cite{Komornikzua}). Symilarly
\begin{equation}
\left\|\begin{array}{cc}
(\Delta v_l, m\cdot \nabla v_l)=(n-2)|| v_l||^2 - \displaystyle \int_{\Gamma}(m\cdot \nu)|\nabla v_l|^2d\Gamma + 2\displaystyle \int_{\Gamma}\frac{\partial u_l}{\partial \nu} (m\cdot \nabla u_l)d\Gamma\\
2(v'_l,m\cdot \nabla v'_l)=-n|v'_l|^2 +\displaystyle \int_{\Gamma}(m\cdot \nu)(v')^{2}d\Gamma. \label{D4}
\end{array}
\right.
\end{equation}

\subsection{Proof of Theorem \ref{teo2}}

We will prove the inequality (\ref{a11}) of theorem \ref{teo2} for solutions $u_l,\;v_l$ given by theorem \ref{teo1} with $h_1(x,s)=(m(x)\cdot \nu(x))p_1(s)$ and $h_2(x,s)=(m(x)\cdot \nu)p_2(s)$. The result follows by taking the infimum limit on both sides of the obtained inequality and using convergence (\ref{a41}).

In order to facilitate the writing, we will omit the sub-index l of the diverse expressions.

Introduce the notation
\begin{eqnarray}
&&\displaystyle E(t)=\frac{1}{2}\Bigl(|u'(t)|^2+\frac{\alpha_1}{\alpha_2}|\sigma'(t)|^2+\mu(t)||u(t)||^2 + \frac{\alpha_1}{\alpha_2}||v(t)||^2\Bigr)\;\;\;t\geq 0\label{D5}
\end{eqnarray}

By similar computations made to obtain (\ref{a15}), we deduce from $(\ref{a37})_1$ and $(\ref{a37})_2$ that after multiplying  $\displaystyle\frac{\alpha_1}{\alpha_2}$ by (\ref{a37}) we obtain
\begin{equation}
\begin{array}{cc}
\displaystyle \frac{dE}{dt}=\frac{u'}{2}||u||^2-\mu \displaystyle \int_{\Gamma}(m, \nu)p_2(v')v'd\Gamma -\alpha_1\displaystyle\Bigl(\sum_{i=1}^{n}\frac{\partial v}{\partial x_i },u'  \Bigr)- \\ \displaystyle \frac{\alpha_1}{\alpha_2} \int_{\Gamma_1}(m, \nu)p_2(v')v'd\Gamma  -\displaystyle \frac{\alpha_1}{\alpha_2} \int_{\Gamma_1} \sigma u v'd\Gamma  +\alpha_1\displaystyle\Bigl(\sum_{i=1}^{n}\frac{\partial u}{\partial x_i },v'  \Bigr)\label{D6}
\end{array}
\end{equation}

By Gauss Theorem, we have
$$\alpha_1\displaystyle\Bigl(\sum_{i=1}^{n}\frac{\partial v'}{\partial x_i },u  \Bigr)=\alpha_1\displaystyle\int{\Gamma_1}uv'\Bigl(\sum_{i=1}^n \nu_i\Bigr)d\Gamma_1 -\alpha_1\displaystyle\Bigl(\sum_{i=1}^{n}\frac{\partial u}{\partial x_i },v'  \Bigr) $$
that implies
 \begin{eqnarray*}
\displaystyle\frac{d}{dt}\Bigl(\alpha_1\sum_{i=1}^n \frac{\partial v}{\partial x_i},u\Bigr)= \alpha_1\displaystyle\int_{\Gamma_1}uv'\Bigl(\sum_{i=1}^n \nu_i\Bigr)d\Gamma_1 -\alpha_1\displaystyle\Bigl(\sum_{i=1}^{n}\frac{\partial u}{\partial x_i },v'  \Bigr)+\alpha_1\displaystyle\Bigl(\sum_{i=1}^{n}\frac{\partial v}{\partial x_i },u'  \Bigr).
\end{eqnarray*}

Therefore
\begin{eqnarray*}
-\alpha_1\displaystyle\Bigl(\sum_{i=1}^{n}\frac{\partial v}{\partial x_i },u'\Bigr)= - \alpha_1\displaystyle\Bigl(\sum_{i=1}^{n}\frac{\partial u}{\partial x_i },v'\Bigr)+\alpha_1\displaystyle\int_{\Gamma_1}uv'\Bigl(\sum_{i=1}^n \nu_i\Bigr)d\Gamma_1 - \displaystyle\frac{d}{dt}\Bigl(\alpha_1\sum_{i=1}^n \frac{\partial v}{\partial x_i},u\Bigr)
\end{eqnarray*}

Combining this equality with (\ref{D6})and canceling similar terms with opposite signs, we obtain
\begin{eqnarray*}
\displaystyle \frac{dE}{dt}+\frac{d}{dt}\Bigl(\alpha_1\displaystyle\sum_{i=1}^{n}\frac{\partial v}{\partial x_i},u\Bigr)=\frac{u'}{2}||u||^2-\mu \displaystyle \int_{\Gamma_1}(m,\nu)p_1(u')u'\Gamma_1-\\
 \frac{\alpha_1}{\alpha_1}\displaystyle \int_{\Gamma_1}(m,\nu)p_2(v')v'd\Gamma_1-\displaystyle\int_{\Gamma_1}\Bigl[\frac{\alpha_1}{\alpha_2}\sigma -\alpha_1\Bigl(\displaystyle \sum_{i=1}^{n}\nu_i    \Bigr) \Bigr]u v'd\Gamma_1
\end{eqnarray*}

Then the hypothesis (\ref{a6}) implies that
\begin{eqnarray*}
\displaystyle \frac{dE}{dt}+\frac{d}{dt}\Bigl(\alpha_1\displaystyle\sum_{i=1}^{n}\frac{\partial v}{\partial x_i},u\Bigr)=\frac{u'}{2}||u||^2-\mu \displaystyle \int_{\Gamma_1}(m,\nu)p_1(u')u'\Gamma_1-\\
 \frac{\alpha_1}{\alpha_1}\displaystyle \int_{\Gamma_1}(m,\nu)p_2(v')v'd\Gamma_1\leq 0
\end{eqnarray*}

Using the notation
\begin{equation}
F(t)=\alpha_1\Bigl(\displaystyle \sum_{i=1}^n\frac{\partial v}{\partial x_i}(t),u(t)\Bigr),\;t\geq 0 \label{D7}
\end{equation}
we obtain
\begin{equation}
\displaystyle \frac{d}{dt}\Bigl( E +F\Bigr)=\frac{u'}{2}||u||^2-\mu \displaystyle \int_{\Gamma_1}(m,\nu)p_1(u')u'\Gamma_1-\frac{\alpha_1}{\alpha_1}\displaystyle \int_{\Gamma_1}(m,\nu)p_2(v')v'd\Gamma_1\label{D8}.
\end{equation}

The above equality provide bounded solutions on $[0, \infty]$. In order to obtain the decay of solutions, we introduce the functional
\begin{equation}
G(t)= (n-1)(u',u) + (n-1)(v',v) +2(u',m.\nabla u)+ 2(v', m.\nabla v). \label{D9}
\end{equation}

\subsection{Boundedness of F an G}

We have
\begin{equation}
|F|\leq 2\Bigl(\displaystyle \frac{\alpha_1.\alpha_2}{\mu_0}n \Bigr)^{\frac{1}{2}}ME \label{D10}
\end{equation}
where $M$ is the constant introduced in (\ref{u}). Then we obtain
$$\begin{array}{l}
{\rm \bullet} \;\;\;{\displaystyle \Bigl|(n-1)(u',u)|\leq 2(n-1)\displaystyle\frac{M}{\mu_0^{\frac{1}{2}}}E}\\[5pt]
{\rm \bullet}\;\;\; {\Bigl|(n-1)(v',v)\Bigr|\leq 2(n-1)M\frac{\alpha_1}{\alpha_2}E}\\[5pt]
{\rm \bullet}\;\;\; {\displaystyle \Bigl|2(u', m.\nabla u) \Bigr|\leq 4\frac{R(x_0)}{\mu_0\frac{1}{2}}E\;\;where\;R(x^0)\;was\;introduced\; in\; (\ref{v}) }\\[5pt]
{\rm \bullet}\;\;\; {\displaystyle \Bigl|2(v',m.\nabla v) \Bigr| \leq 4 R(x_0)\frac{\alpha_1}{\alpha_2}E}
\end{array}$$

Thus
\begin{equation}
\bigl|G\bigr|\leq AE  \label{D11}
\end{equation}
Where the constant $A$ was introduced in (\ref{a1}).
From (\ref{D10}), (\ref{D11}) and $\varepsilon $, it follows that
$$ |F+\varepsilon G|\leq \Bigl[2\bigl(\frac{\alpha_1 \alpha_2 n}{\mu_0}\bigr)^{\frac{1}{2}} M + \varepsilon A \Bigr]E $$

For the particular $\alpha_1\;\alpha_2$ satisfying hypothesis (\ref{a5}), we have
$$|F+ \varepsilon_1 G|\leq \frac{1}{2}E,\;\; 0\leq\varepsilon_1\leq \frac{1}{4s} $$
So
\begin{equation}
\frac{1}{2}E(t)\leq E(t) + F(t)+ \varepsilon_1 G(t)\leq \frac{2}{3}E(t),\; t\geq 0.   \label{D12}
\end{equation}

\subsection{Boundedness of $\frac{dG}{dt}$}

We have
\begin{equation}
\begin{array}{cc}
\displaystyle \frac{dG}{dt}=(n-1)(u^{''},u)+(n-1)|u'|^2+(n-1)(v^{''},v)+(n-1)|v'|^2+\\ 2(u^{''},m.\nabla u)+ 2(u',m.\nabla u')+ 2(v^{''},m.\nabla v)+2(v',m.\nabla v') =\;\;\;\;\;\\
I_1+(n-1)|u'|^2 + I_2+ (n-1)|v'|^2 +I_3 + I_4 + I_5 + I_6\;\;\;\;\;\;\;\;\;\; \label{D13}
\end{array}
\end{equation}

$\bullet$ By equation $(\ref{a37})_1$, we find $$I_1=(n-1)(\mu \Delta u,u) - \alpha_1(n-1)(\displaystyle\sum_{i=1}^n\frac{\partial v}{\partial x_i},u)$$ and by equation $(\ref{a40})_1$, we also find $$I_1=-(n-1)\mu||u||^2-(n-1)\mu \displaystyle \int_{\Gamma_1}(m.\nu)h_1(u')ud\gamma_1-(n-1)\alpha_1\Bigl(\displaystyle \sum_{i=1}^n\frac{\partial v}{\partial x_i}, u\Bigr) $$

$\bullet$ In a similar way, by $(\ref{a37})_2$ and $(\ref{a40})_2$, we derive,

\begin{eqnarray*}
 &&\displaystyle I_2=-(n-1)||v||^2-(n-1)\displaystyle \int_{\Gamma_1}(m.\nu)h_2(v')vd\Gamma_1- (n-1)\displaystyle \int_{\Gamma_1}(\sigma u)vd\Gamma_1+ \\
&&\displaystyle (n-1)\alpha_2\Bigl(\displaystyle \sum_{i=1}^n\frac{\partial u}{\partial x_i}, v\Bigr)
\end{eqnarray*}

$\bullet$ By equation $(\ref{a37})_1$ and identity $(\ref{D3})$, we get

\begin{eqnarray*}
 &&\displaystyle I_3=\mu(n-1)||u||^2-(n-1)\displaystyle \int_{\Gamma_1}(m.\nu)|\nabla u|^2d\Gamma_1+ \displaystyle 2\mu\int_{\Gamma_1}\frac{\partial u}{\partial x_i}(m.\nabla u)d\Gamma_1-\\
&&\displaystyle 2\alpha_1\Bigl(\displaystyle \sum_{i=1}^n\frac{\partial v}{\partial x_i}, m.\nabla u\Bigr)
\end{eqnarray*}

$I_4=-|u'|^2+\displaystyle\int_{\Gamma}(m.\nu)u^{\frac{1}{2}}d\Gamma_1.$

$\bullet$ In a similar way, by $(\ref{a37})_2$ and (\ref{D3}), we find

\begin{eqnarray*}
 &&\displaystyle I_5=\mu(n-1)||v||^2-\displaystyle \int_{\Gamma_1}(m.\nu)|\nabla v|^2 d\Gamma_1+ \displaystyle 2\int_{\Gamma_1}\frac{\partial v}{\partial x_i}(m.\nabla v)d\Gamma_1+\\
&&\displaystyle 2\alpha_2\Bigl(\displaystyle \sum_{i=1}^n\frac{\partial u}{\partial x_i}, m.\nabla v\Bigr)
\end{eqnarray*}

$I_6= -n|v'|^2 +\displaystyle \int_{\Gamma_1}(m.\nu)(v')^2d\Gamma_1. $

Taking into account the last four equalities in (\ref{D13}) and canceling  the terms with opposite signs, we have
\begin{equation}
\begin{array}{cc}
\displaystyle \frac{dG}{dt}= -|u'|^2-|v'|^2-\mu ||u||^2-||v||^2-(n-1)\mu\displaystyle \int_{\Gamma}(m.\mu)h_1(u')ud\Gamma_1 -\\(n-1)\alpha_1\Bigl(\displaystyle\sum_{i=1}^n \frac{\partial v}{\partial x_i}, u\Bigr)-(n-1)\displaystyle\int_{\Gamma_1}(m.\nu)h_2(v')vd\Gamma_1 -(n-1)\displaystyle \int_{\Gamma_1}\sigma uvd\Gamma_1+ \\
(n-1)\alpha_2\Bigl(\displaystyle\sum_{i=1}^n \frac{\partial u}{\partial x_i}, v\Bigr) - \mu \displaystyle \int_{\Gamma_1}(m.\nu)|\nabla v|^2d\Gamma + 2\displaystyle \mu \int{\Gamma}\frac{\partial u}{\partial \nu}(m.\nabla u)d\Gamma- \\
2\alpha_1\Bigl(\displaystyle\sum_{i=1}^n \frac{\partial v}{\partial x_i}, m.\nabla u \Bigr)- \displaystyle\int_{\Gamma_1}(m.\nu)|\nabla v|^2d\Gamma_1 + 2 \displaystyle\int_{\Gamma_1}\frac{\partial v}{\partial x_i}(m.\nabla v) d\Gamma_1 + \\
2\alpha_2 \Bigl(\displaystyle\sum_{i=1}^n \frac{\partial u}{\partial x_i}, m.\nabla v \Bigl)+ \displaystyle \int_{\Gamma_1}(m.\nu)(u')^2d\Gamma_1 + \displaystyle \int_{\Gamma_1}(m.\nu)(v')^2d\Gamma_1 = \\
-|u'|^2 - |v'|^2 -\mu ||u||^2-||v||^2+\displaystyle \sum_{k=1}^{11}J_{k} + \displaystyle \int_{\Gamma_1}(m.\nu)(u')^2d\Gamma +\\
\displaystyle \int_{\Gamma_1}(m.\nu)(v')^2d\Gamma \;\;\;\;\;\;\;\;\;\;\;\;\;\;\;\;\;\;\;\;\;\;\;\;\;\;\;\;\;\;\;\;\;\;\;\;\;\;\;\;\;\;\;\;\;\;\;\;\;\;\;\;\;\;\;\;\;\;\;\;\;\
\label{D14}
\end{array}
\end{equation}

we also have,
$$ \bullet \;\; |J_1|\leq 4(n-1)^2\mu(0)R(x^0)L_1^2N^2\displaystyle\int_{\Gamma_1}(m.\nu)(u')^2d\Gamma_1 +\frac{1}{16}\mu||u||^2 $$
where the constants $L_i\; (i=1,\;2)$ were introduced in hypothesis (\ref{y})
$$ \bullet \;\; |J_2|\leq 4(n-1)^2nM^2\displaystyle\frac{\alpha_1^2}{\mu_0}(\mu||u||^2)+\frac{1}{16}||v||^2 \;\;\;\;\;\;\;\;\;\;\;\;\;\;\;\;\;\;\;\;\;\;\;\;\;\;\;\;$$

$$ \bullet \;\; |J_3|\leq 4(n-1)^2R(x^0)L_2^2 N^2\displaystyle \int_{\Gamma_1}(m.\nu)(v')^2d\Gamma +\frac{1}{16}||v||^2 \;\;\;\;\;\;\;\;\;\;$$

The hypothesis (\ref{a6}) provides

$$ \bullet \;\; |J_4|\leq 4(n-1)^2\Bigl|\Bigl|\displaystyle\sum_{i=1}^n \nu \Bigr|\Bigr|_{L^2{\infty}(\Gamma_1)} \frac{ N^4}{ \mu_0}\alpha^2_2(\mu ||u||^2)+\frac{1}{16}||v||^2 \;\;\;\;\;\;\;\;\;$$

$$ \bullet \;\; |J_5|\leq 4(n-1)^2nM^2 \frac{ \alpha_2}{\mu_0}(\mu ||u||^2)+\frac{1}{16}||v||^2 \;\;\;\;\;\;\;\;\;\;\;\;\;\;\;\;\;\;\;\;\;\;\;\;\;\;\;\;$$

 Observing that $|\nabla u|^2 = (\frac{\partial u}{\partial \nu})^2 $ in $ \Gamma_0 $, we find

\begin{equation}
\bullet J_6 = - \mu \displaystyle \int_{\Gamma_0}(m.\nu)(\frac{\partial u}{\partial \nu})^2d\Gamma_0-\mu \displaystyle \int_{\Gamma_0}(m.\nu)|\nabla u|^2d\Gamma_1\;\;\;\;\;\;\;\;\;\;\;\;\;\; \label{D15}
\end{equation}

noting that $\nabla u = \nu.\frac{\partial u}{\partial x_i}$ on $\Gamma_0$, we obtain

$$J_7 = 2\mu \displaystyle\int_{\Gamma_0}(m.\nu)(\frac{\partial u}{\partial \nu})^2 d\Gamma_0+2\mu \displaystyle\int_{\Gamma_0}\frac{\partial u}{\partial \nu}(m.\nabla u) d\Gamma_1  $$

On the other hand, using equations $(\ref{a37})_1$ and the hypothesis (\ref{v}), (\ref{a10}), we get
\begin{eqnarray*}
 && \Bigl|2\mu \displaystyle\int_{\Gamma_0}\frac{\partial u}{\partial \nu}(m.\nabla u) d\Gamma_1\Bigr|\leq\mu(0)R^2(x^0)\displaystyle\int_{\Gamma_1}\frac{1}{m.\nu}(\frac{\partial u}{\partial x_i})^2 d\Gamma_1 +  \displaystyle\int_{\Gamma_0}(m.\nu)|\nabla u|^2d\Gamma_1 \\
 &&\mu(0)R^2(x^0)L_1^2\displaystyle\int_{\Gamma_1}(m.\nu)(u')^2d\Gamma_1+\mu\displaystyle\int_{\Gamma_1}(m.\nu)|\nabla u|^2d\Gamma_1
 \end{eqnarray*}

So
\begin{equation}
\begin{array}{cc}
J_7\leq 2\mu \displaystyle\int_{\Gamma_0}(m.\nu)(\frac{\partial u}{\partial \nu})^2d\Gamma_0 +\mu(0)R^2(x^0)L_1^2\displaystyle\int_{\Gamma_1}(m.\nu)(u')^2d\Gamma_1+\\
\mu\displaystyle\int_{\Gamma_1}(m.\nu)|\nabla u|^2d\Gamma_1 \;\;\;\;\;\;\;\;\;\;\;\;\;\;\;\;\;\;\;\;\;\;\;\;\;\;\;\;\;\;\;\;\;\;\;\;\;\;\;\;\;\;\;\;\;\; \;\;\;\;\;\;\;\;\label{D16}
\end{array}
\end{equation}

Therefore, after adding (\ref{D15}) and (\ref{D16}), reducing similar terms, canceling similar terms with opposite signs  and noting that $m.\nu \leq 0$ on $\Gamma_0$, we obtain
$$\;\; J_6 +J_7 \leq \mu(0) R^2(x^0)L_1^2\displaystyle\int_{\Gamma_1}(m.\nu)(u')^2d\Gamma_1  $$

$\bullet\;\;\;|J_8|\leq 16R^2(x^0)\frac{n}{\mu_0}\alpha_1^2 (\mu||u||^2)+\frac{1}{16}||v||^2 $

$\bullet$ In a similar way as in (\ref{D16}), we find

\begin{equation}
\begin{array}{cc}
J_9 +J_{10} \leq R^2(x^0)\displaystyle\int_{\Gamma_1}\frac{1}{m.\nu}[-(m\nu)h_2(v')-\sigma u]^2d\Gamma_1 \leq \\
2R^2(x^0)L_2^2\displaystyle\int_{\Gamma_1}(m.\nu)(v')^2d\Gamma_1 +2\alpha^2_2\displaystyle\int_{\Gamma_1}\frac{1}{m.\nu}\Bigl(\displaystyle \sum_{i=1}^n\nu_i\Bigr)^2u^2d\Gamma_1, \label{D17}
\end{array}
\end{equation}

we obtain
$$2\alpha^2_2\displaystyle\int_{\Gamma_1}\frac{1}{m.\nu}\Bigl(\displaystyle \sum_{i=1}^n\nu_i\Bigr)^2u^2d\Gamma_1\leq 2||\displaystyle\sum_{i=1}^n\nu_i||_{L^{\infty}(\Gamma_1)}^2\frac{N^2}{\tau_0\mu_0}\alpha_2(\mu||u||^2)  $$

where the constant $\tau_0$ was introduced in (\ref{v}). The preceding two inequalities provides

$$J_9 + J_{10}\leq 2R^2(x_0)L_2^2\displaystyle\int_{\Gamma_1}(m.\nu)(v')^2d\Gamma_1+ 2||\displaystyle \sum_{i=1}^n\nu_i||_{L^{\infty}(\Gamma_1)}^2\frac{N^2}{\tau_0\mu_0}(\mu||u||^2)$$

$\bullet\;\; |J_{11}|\leq 16R^2(x^0)\frac{n}{\mu_0}\alpha_2^2(\mu||u||^2)+\frac{1}{16}||v||^2 $

Taking into account the above boundedness for $J_i,\;(i=1,...,11)$ in (\ref{D14}) and using notations introduced in (\ref{a2})- (\ref{b4}) we obtain,
\begin{eqnarray*}
&&\displaystyle\frac{dG}{dt}\leq -|u'|^2-|v'|^2-\mu||u||^2-||v||^2+p_1\alpha_1^2(\mu||u||^2)+p_2\alpha_2^2(\mu||u||^2)+\\
&&\displaystyle S_1\int_{\Gamma_1}(m.\nu)(u')^2d\Gamma_1+\displaystyle S_2\int_{\Gamma_1}(m.\nu)(v')^2d\Gamma_1+\frac{1}{16}\mu||u||^2+\frac{3}{8}||v||^2
\end{eqnarray*}

This implies
\begin{eqnarray*}
&&\displaystyle\frac{dG}{dt}\leq -|u'|^2-|v'|^2-\mu||u||^2-||v||^2+p_1\alpha_1^2(\mu||u||^2)+(p_1\alpha_1+p_2\alpha_2)(\mu||u||^2)+\\
&&\displaystyle S_1\int_{\Gamma_1}(m.\nu)(u')^2d\Gamma_1+\displaystyle S_2\int_{\Gamma_1}(m.\nu)(v')^2d\Gamma_1
\end{eqnarray*}

The hypothesis (\ref{a5}) provide $ \frac{15}{16}-\Bigl[p_1\alpha_1+p_2\alpha_2^2\bigr]\geq \frac{1}{2}$. Then
\begin{eqnarray*}
&&\displaystyle\frac{dG}{dt}\leq -|u'|^2-|v'|^2-\mu||u||^2-||v||^2+S_1\int_{\Gamma_1}(m.\nu)(u')^2d\Gamma_1+\displaystyle S_2\int_{\Gamma_1}(m.\nu)(v')^2d\Gamma_1
\end{eqnarray*}

We note that $\frac{1}{2}\frac{\alpha_1}{\alpha_2}\leq \frac{1}{2}$ or $-\frac{1}{2}< -\frac{1}{2}\frac{\alpha_1}{\alpha_2}$ for all $\alpha_1>0$ and $\alpha_2>0$. Thus

\begin{eqnarray}
&& \displaystyle \frac{d G}{dt}\leq -E + S_1 \displaystyle\int_{\Gamma_1}(m.\nu)(u')^2d\Gamma_1 +\displaystyle\int_{\Gamma_1}(m.\nu)(v')^2d\Gamma_1\label{D18}
\end{eqnarray}

In the sequel, we conclude the proof of Theorem \ref{teo2}. By (\ref{D8}), (\ref{D17}), hypothesis (\ref{y}) and for $\epsilon >0$, we have

 \begin{eqnarray}
&& \displaystyle \frac{d}{dt}(E+F+\epsilon G)\leq \frac{\mu'}{2}||u||^2-\mu_0b_1\displaystyle  \int_{\Gamma_1}(m.\nu)(u')^2d\gamma_1-\nonumber\\
&& \displaystyle\frac{\alpha_1b_2}{\alpha_2}\displaystyle \int_{\Gamma_1}(m.\nu)(v')^2d\Gamma_1-\epsilon E+\epsilon S_1\displaystyle \int_{\Gamma_1}(m.\nu)(u')^2d\Gamma_1\\ \label{D19}
&& \displaystyle \epsilon S_2\int_{\Gamma_1}(m.\nu)(v')^2d\Gamma_1\nonumber
\end{eqnarray}

Choosing $\epsilon_2>0$ in conditions (\ref{a7}) we find
$$\frac{d}{dt}(E+F+\epsilon_2G)\leq -\epsilon_2E  $$

taking $\eta >0$ in conditions (\ref{a8}) and using (\ref{D12}),we get

$$E(t)+F(t)+\eta G(t)\leq e^{-\frac{2}{3}\eta}\Bigl[E(o)+F(0)+\eta F(0)\Bigr]$$

Then (\ref{D12}) implies that

$$E(t)\leq 3e^{-\frac{2}{3}\eta}E(0),\;\;\forall t\geq 0  $$

the proof is completed.

\acknowledgement{The author M. L. Oliveira acknowledges the support of National
Institute of Science and Technology of Mathematics INCT-Mat and
CAPES and CNPq/Brazil.}


\begin{thebibliography}{99}
  \bibitem{Araruna} Araruna, F. D. and Maciel A. B., {\it Existence and boundary stabilization of the semilinear wave equations}, Nonlinear Anal. 67(2007), 1288 - 1305.

  \bibitem{Araujo} Araujo, J. L. G., Milla Miranda, M. and Medeiros, L. A., {\it Vibrations of beam by torsion or impact}, Math. Comtemp. 36(2009), 29 - 50

 \bibitem{Cavalcanti} Cavalcanti, M. M., Calvalcanti, V. N. D. and Martinez, P.,{\it Existence and decay rate estimates for the wave equation with nonlinear boundary damping and source term}, J. Diff. Eq. 203(2004), 114 -158.

  \bibitem{Clark} Clark, H. R.,San gil Jutuca, L. P. and Milla Miranda,M., {\it On a mixed problem for a linear coupled system with variable coefficients}, Eletronic J. Diff. Equations 4(1998), 1- 20

  \bibitem{Kim} Kim, J. U.. and Renardy, Y.,{\it Boundary control of Timoshenko beam}, SIAM.J Control and Optimization, 25 (1987), 1417 - 1429.

 \bibitem{Komornik} Komornik, V., {\it Exact controllability and stabilization - Multiplier Method}, J. Wiley and Masson,Paris, 1944.

 \bibitem{Komornikzua}Komornik, V. and Zuazua, E.,{\it A direct method for boundary stabilization of the wave equation}, J. Math. Pure et Appl. 69(1990), 33 - 54.

  \bibitem{Laz} Lasiecka, I. and Tataru, D., {\it Uniform boundary stabilization of semilinear wave equation with nonlinear damping}, Diff. Integral Eq. 6(1993), 507 - 533.

  \bibitem{Lions} Lions, J. L.,{\it Quelques M\'ethodes de R\'esolutions des Probl\'emes aux Limites Non-lin\'eares }, Dunod, Paris, 1964.

  \bibitem{LJL} Lions, J. L.,{\it \'Equations aux D\'eriv\'ees Partielles - Interpolation Vol. I}, EDP sciences, Les Ulis, Paris, France 2003, {it Oeuvres choisis de Jacques-Lions}, 2003, view at Math sci Net.

 \bibitem{LM} Louredo, A. T, and Milla Miranda, M., {\it Nonlinear boundary dissipation for a coupled system of Klein-Gordon equations}, Electronic J. Diff. Eq. 120(2010), 1-19.

  \bibitem{Louredo} Louredo, A. T. and Milla Miranda, M {\it Local solutions for a coupled system of Kirchhoff type}, Nonlinear Anal. 74(2011), 7094 - 7110.

  \bibitem{millaaraujo} A. T. Louredo, M. A. F. and Milla Miranda, M., {\it On a nonlinear wave equation with boundary damping}, Math. Methods Appl. Sci., Accepted to publications, DOI: 10.1002/mma. 2885.

   \bibitem{Marcus} Marcus, M. and Miael, V., {\it Every superposition operator mapping one Sobolev space into another is continuous}, J. Funct. anal. 33(1979), 217 - 229.

  \bibitem{Milla} Milla Miranda, M. and Medeiros, L. A., {\it On a boundary value problem for wave equations: Existence-Uniqueness-Asymptotic behavior}, Res. Mat. Apl. Univ. do Chile, 17(1446), 47-73.

  \bibitem{Millajutuca} Milla Miranda, M. and san Gil Jutuca, L. P.,{\it Existence and boundary stabilization of solutions for the Kirchhoff equation}, Commum. Partial Diff. Eq. 24(1999), 1759 - 1800

  \bibitem{Mota} Mota, A.C. P. C.,{\it Existence and Stability for a Nonlinear timoshenko system}, Ph. D. Dissertation, IM - UFRJ, Rio, 2004.

  \bibitem{Strauss} Strauss, W. A., {\it On work solutons of semilinear hyperbolic equations}, An. Acad. Bras. Ci\^enc. 42(1970), 645 - 651.
  
   \bibitem{Timoshenko} Timoshenko, S. and Woinowsky Krieger, S., {\it Theory of Plates and Shells}, MacGraw-Hill, New York, (1959).

  \bibitem{Tucsnak} Tucsnak,S., {\it On an initial and boundary value problem for the nonlinear timoshenko Beam}, An. Acad. Brs. Ci\^enc. 63(1991), 115 - 125.
  
   \bibitem{Vitilaro} Vitillaro, E., {\it Global existence for the wave equation with nonlinear boundary damping and source terms}, J, Diff. Eq. 186(2002), 259- 298.

   \bibitem{Zuazua} Zuazua, E., {\it Uniform stabilization of the wave equation by nonlinear feedback}, SIAM J. control Optim. 28(1990), 466 - 478.


\end{thebibliography}
\end{document}